\documentclass[10pt]{amsart}
\usepackage{amsfonts, mathrsfs}
\usepackage{ifthen}
\usepackage{amsthm}
\usepackage{amsmath}
\usepackage{graphicx}
\usepackage{amscd,amssymb,amsthm}
\usepackage{graphicx}
\usepackage{epstopdf}
\usepackage{hyperref}
\usepackage{algorithm2e}
\usepackage{mathtools}
\usepackage{color}
\newcounter{minutes}
\divide\time by 60
\newcounter{hours}
\multiply\time by 60 \addtocounter{minutes}{-\time}

\newtheorem{lemma}{Lemma}
\newtheorem{theorem}{Theorem}

\newtheorem{remark}{Remark}

\newcommand{\real}{\operatorname{Re}}

\keywords{Bessel functions; convex functions; radius of convexity; radius of uniform convexity; Rayleigh sums;
asymptotic expansion; Euler-Rayleigh inequalities; Weierstrassian factorization; asymptotic inversion; ordinary potential polynomials.}
\subjclass[2010]{30C45, 33C10.}

\begin{document}
\title[Radii of convexity and uniform convexity]{Bounds and asymptotic expansions for the radii of convexity and uniform convexity of normalized Bessel functions}
\allowdisplaybreaks
%=======================================================================================================================================================
\author[\'A. Baricz]{\'Arp\'ad Baricz}
\address{Department of Economics,  Babe\c{s}-Bolyai University, Cluj-Napoca 400591, Romania}
\address{Institute of Applied Mathematics, \'Obuda University, 1034 Budapest, Hungary}
\email{bariczocsi@yahoo.com}

\author[P. Kumar]{Pranav Kumar}
\address{Department of Mathematics,
Indian Institute of Technology Madras, Chennai 600036, India}
\email{pranavarajchauhan@gmail.com}

\author[S. Singh]{Sanjeev Singh}
\address{Department of Mathematics,
Indian Institute of Technology Indore, Indore 453552, India}
\email{snjvsngh@iiti.ac.in}
\newcommand{\ec}{\operatorname{erfc}}

\def\thefootnote{}
\footnotetext{ \texttt{File:~\jobname .tex,
         printed: \number\year-0\number\month-\number\day,
          \thehours.\ifnum\theminutes<10{0}\fi\theminutes}
} \makeatletter\def\thefootnote{\@arabic\c@footnote}\makeatother

\maketitle

\begin{abstract}
This paper explores the asymptotic behavior of the radii of convexity and uniform convexity for normalized Bessel functions with respect to large order. We provide detailed asymptotic expansions for these radii and establish recurrence relations for the associated coefficients. Additionally, we derive generalized bounds for the radii of convexity and uniform convexity by applying the Euler-Rayleigh inequality and potential polynomials. The asymptotic inversion method and Rayleigh sums are the main tools used in the proofs.
\end{abstract}
%====================================================================================================================================
\section{Introduction}

Bessel functions, long regarded as fundamental in the realm of classical special functions, hold a pivotal role in mathematical analysis, physics, and engineering. Their geometric properties, deeply rooted in complex function theory, became a prominent area of exploration during the 1960s, thanks to the pioneering work of scholars such as Brown, Hayden, Kreyszig, Merkes, Scott, Robertson, and Wilf \cite{Br60, Br62, Br82, HM64, KT60, MRS62, Ro54, Wi62}. In recent years, this investigation has intensified, with a focus on characteristics like univalence, starlikeness, and convexity \cite{ABY17, BKS14, BOS16, BS14, DS17, Sz15}. For normalized Bessel functions of the first kind, significant progress has been made in determining radii and orders of starlikeness and convexity.

Baricz and Szász \cite{BS14} used
Mittag-Leffler expansions
for quotients of Bessel functions and the fact that smallest positive zeros of some Dini functions are less
than the first positive zero of the Bessel functions of the first kind to deduce the radius of convexity for normalized Bessel functions. Similarly, Akta\c{s} et al. \cite{ABO18} utilized Euler-Rayleigh inequalities to establish bounds for the radii of starlikeness and convexity for normalized Bessel, Struve, and Lommel functions. More recently, Baricz and Nemes \cite{BN21} introduced a systematic asymptotic expansion for the radius of starlikeness of normalized Bessel functions, offering deeper insights into their geometric behavior.

Motivated by these advancements, this paper extends the study of normalized Bessel functions of the first kind by deriving comprehensive asymptotic expansions for the radii of convexity and uniform convexity with respect to large orders. Central to our methodology are the asymptotic properties of Rayleigh sums and the Laurent series expansions for the positive zeros of Bessel functions at infinity. We express the coefficients of these expansions using ordinary potential polynomials and provide recurrence relations to support the computations.

Beyond their geometric properties, the zeros of Bessel functions have numerous applications, including wave propagation, scattering theory, and quantum mechanics \cite{DYL6,ELR93,FS08,LZ07,Pa72}. These diverse applications have inspired extensive research on the asymptotic expansion of the zeros of Bessel functions. Applications and recent contributions  in this area can be found in \cite{QW99, Du24} and references therein. The radii of convexity and uniform convexity of Bessel functions are determined by equations involving Bessel functions and their derivatives. Their asymptotic expansions are particularly useful in approximating the zeros of these functions for large values of $\nu$ by truncating the series to a finite number of terms. Furthermore, in cases where such functions are used to approximate others, these asymptotic expansions can significantly aid in approximating the zeros of other functions.

Our approach not only complements the work of Baricz and Nemes \cite{BN21} but also offers generalized bounds for the radii of convexity and uniform convexity in terms of potential polynomials. Additionally, graphical representations of the approximate radii for large orders substantiate our theoretical findings, illustrating that for fixed but large orders, the radius of uniform convexity remains smaller than the radius of convexity for a given normalized Bessel function.These results contribute to a richer understanding of the complex geometric properties of Bessel functions and their broader implications.

Before stating our results, we introduce some necessary notations and definitions. Consider the set $\mathbb{D}_r$ defined as $\mathbb{D}_r = \{z \in \mathbb{C}: |z| <r\}$, where $r > 0$. Let $f : \mathbb{D}_r \to \mathbb{C}$ be a normalized univalent or one-to-one function that satisfies the conditions $f(0) = 0$ and  $f '(0) = 1$. In other words, $f$ takes the form $f(z) = z + a_2z^2 + a_3z^3 +\ldots$, where the coefficients $a_2, a_3, \ldots$ are real or complex numbers. The radius of univalence of the function $f$ is the largest radius $r$ for which $f$ maps univalently the open disk $\mathbb{D}_r$ into some domain in the complex plane. Similarly, the radius of convexity of the function $f$ is the largest radius $r$ for which $f$ maps $\mathbb{D}_r$ into a convex domain. It is worth to mention that the class of normalized convex functions (with respect to the origin) is a subclass of univalent functions. Consequently, the radius of univalence of $f$ is greater than or equal to the radius of convexity of the same function $f$. Considering the analytic characterization of convex functions, the radius of convexity is determined by
$$r^c(f)=\text{sup}\left\{r\in \left(0,\infty\right)\left|\real\left(1+\frac{zf''(z)}{f'(z)}\right)>0\quad \mbox{for all}\quad z\in \mathbb{D}_r\right.\right\}.$$

The concept of uniform convexity was introduced by Goodman \cite{Go91}. A function $f$ is said to be uniformly convex in $\mathbb{D}_r$ if $f$ is a convex function and has the property that every circular arc $\gamma$ contained in $\mathbb{D}_r$, with center $\xi$, the arc $f(\gamma)$ is convex. Analytically (see \cite{Ro93} or \cite[Theorem 2.1]{DS17}), the function $f(z)=z+z + a_2z^2 + a_3z^3 +\ldots$ in the disk $\mathbb{D}_r$ is uniformly convex if and only if
\begin{align*}
	\real\left(1+\frac{zf^{\prime\prime}(z)}{f^\prime(z)}\right)>\left|\frac{zf^{\prime\prime}(z)}{f^\prime(z)}\right|\quad \mbox{for all}\quad z\in \mathbb{D}_r.
\end{align*}
The radius of uniform convexity is defined by
\begin{align*}
	r^{uc}(f)=\sup\left\{r\in (0,r^c(f))\left|\real \left(1+\frac{zf^{\prime\prime}(z)}{f^\prime(z)}\right)>\left|\frac{zf^{\prime\prime}(z)}{f^\prime(z)}\right|\quad \mbox{for all}\quad z\in \mathbb{D}_r\right.\right\}.
\end{align*}

Now, we turn our attention to the Bessel function of the first kind of order $\nu,$ which is defined by \cite[p. 217]{OLBC10}
\begin{equation}\label{Bessel_series}
	J_\nu(z)=\sum_{n\geq0}^{}\frac{(-1)^nz^{2n+\nu}}{2^{2n+\nu}n!\Gamma(n+\nu+1)},
\end{equation}
and its derivative
\begin{equation}\label{der_Bessel_series}
	J_\nu^{\prime}(z)=\sum_{n\geq0}^{}\frac{(-1)^n(2n+\nu)z^{2n+\nu-1}}{2^{2n+\nu}n!\Gamma(n+\nu+1)},
\end{equation}
respectively.

In this paper, our attention is directed toward the following two normalized forms

\begin{equation}\label{unif_g_nu_def}
g_\nu(z)=2^\nu\Gamma(\nu+1)z^{1-\nu}J_\nu(z)=z-\frac{1}{4(\nu+1)}z^3+\frac{1}{32(\nu+1)(\nu+2)}z^5-\ldots,
\end{equation}
\begin{equation}\label{unif_h_nu_def}
	h_\nu(z)=2^\nu\Gamma(\nu+1)z^{1-\frac{\nu}{2}}J_\nu(\sqrt{z})=z-\frac{1}{4(\nu+1)}z^2+\ldots,
\end{equation}
where $\nu>-1$. It is important to mention that $g_\nu(z)$ and $h_\nu(z)$ remain well-defined even when $\nu< -1$ and not a negative integer. However, in our paper, the condition $\nu > -1$ is of great significance, as it ensures that the zeros of the Bessel function $J_\nu$ are all real, as stated in \cite[p. 482]{Wa44}. The reality of these zeros is a crucial factor in our paper because all results from \cite{ABO18} and \cite{BS14} concerning the radii of convexity, which we will be using, rely on this condition.

Before we start to present our main results, we first introduce some Dini functions which play an important role in the proofs. Let us consider the Dini function $ d_\nu:\Omega\subseteq \mathbb{C}\to\mathbb{C}$ and $e_\nu:\Omega\subseteq \mathbb{C}\to\mathbb{C},$ defined by
\begin{equation}\label{dini_d}
d_\nu(z)=(1-\nu)J_\nu(z)+zJ_\nu^\prime(z)
\end{equation}
and
\begin{equation}\label{dini_e}
e_\nu(z)=(2-\nu)J_\nu(z)+zJ_\nu^\prime(z).
\end{equation}
We define the Rayleigh functions or Rayleigh sums associated with zeros of $d_\nu(z)$ and $e_\nu(z)$ as $\eta_k(\nu)$ and $\theta_k(\nu)$, respectively, with the following formulation
\begin{equation}\label{eta_def}
\eta_k(\nu)=\sum_{n=1}^{\infty}\frac{1}{\alpha_{\nu,n}^{2k}}
\end{equation}
and
\begin{equation}\label{theta_def}
	\theta_k(\nu)=\sum_{n=1}^{\infty}\frac{1}{\beta_{\nu,n}^{2k}}.
\end{equation}
Here, $k$ represents any positive integer, and $\nu$ is a real number with $\nu > -1$. The symbols $\alpha_{\nu,n}$ and $\beta_{\nu,n}$ denote the $n$th positive zero of the Dini functions $d_\nu(z)$ and $e_\nu(z)$, respectively. It is worth noting that owing to their significance in problems associated with Bessel functions, these Rayleigh sums may hold an independent interest. Furthermore, throughout this paper, unless explicitly specified otherwise, vacant summations are considered to be equivalent to zero. Additionally, $\mathbb{N}$ is the set of all positive integers and $\mathbb{N}_0= \mathbb{N}\cup\{0\}$.

The Laguerre-P\'{o}lya class of entire functions (denoted by $\mathcal{LP}$) serves a pivotal role in deriving bounds for the radius of uniform convexity. A real entire function $\phi$ belongs to the class $\mathcal{LP}$ if it can be represented in the form
\begin{align*}
	\phi(z)=cz^d e^{-\alpha z^2+\beta z}\prod_{n\geq1}^{}\left(1-\frac{z}{z_n}\right)e^{z/z_n},
\end{align*}
where $c, \beta,z_n\in \mathbb{R}$, $\alpha\geq 0$, $d\in \mathbb{N}_0$ and $\sum_{n\geq1}^{}z_n^{-2}<\infty$. An important property of this class is that it is closed under differentiation, meaning that if $\phi \in \mathcal{LP}$, then $\phi^{(m)} \in \mathcal{LP}$ for all non-negative integers $m$. For a deeper understanding of the $\mathcal{LP}$ class, readers are referred to \cite[p. 703]{DC09} and references therein.

The paper is organized as follows: Section \ref{sectionA} outlines the main results and lemmas, presenting some asymptotic expansions for the radii of convexity and uniform convexity for two kinds of normalized Bessel functions of the first kind. Section \ref{proofs} contains the proofs of the main theorems and the lemmas.

\section{Preliminary and main results}\label{sectionA}
\setcounter{equation}{0}
\subsection{Asymptotic expansions for the radii of convexity of normalized Bessel functions}

We start with some basic results concerning the above mentioned Rayleigh sums.

\begin{lemma}\label{lemma_eta}
For any positive integer $k$ and positive real $\nu>k$, the Rayleigh sum in $\eqref{eta_def}$ has the convergent Laurent expansion
\begin{equation}\label{etaseries}
\eta_k(\nu)=\frac{1}{\nu^k}\sum_{n=0}^{\infty}\frac{\eta_n^{(k)}}{\nu^n},
\end{equation}
where for any fixed non negative integer $n$, the coefficients $\eta_n^{(k)}$ can be evaluated by the recurrence relation
\begin{align*}
	\eta_n^{(k)}=-ka_n^{(k)}-\sum_{m=0}^{n}\sum_{i=1}^{k-1}a_m^{(i)}\eta_{n-m}^{(k-i)}
\end{align*}
and $a_n^{(k)}$ is given by
$$a_{n}^{(k)}=\frac{\left(-1\right)^k\left(2k+1\right)}{2^{2k}k!}\sum_{k_{k-1}=0}^{n}\ldots\sum_{k_2=0}^{k_3}\sum_{k_1=0}^{k_2}\left(-1\right)^{k_n}\left(-2\right)^{k_2-k_1}\ldots\left(-k\right)^{n-k_{k-1}},\qquad n\in \mathbb{N}_0.$$
\end{lemma}

\begin{lemma}\label{lemma_theta}
For any positive integer $k$ and positive real $\nu>k$, the Rayleigh sum in $\eqref{theta_def}$ has the convergent Laurent expansion
\begin{equation}\label{thetaseries}
\theta_k(\nu)=\frac{1}{\nu^k}\sum_{n=0}^{\infty}\frac{\theta_n^{(k)}}{\nu^n},
\end{equation}
where for any fixed non negative integer $n$, the coefficients $\theta_n^{(k)}$ can be evaluated by the recurrence relation
\begin{align*}
\theta_n^{(k)}=-kb_n^{(k)}-\sum_{m=0}^{n}\sum_{i=1}^{k-1}b_m^{(i)}\theta_{n-m}^{(k-i)}
\end{align*}
and $b_n^{(k)}$ is given by
$$b_{n}^{(k)}=\frac{\left(-1\right)^k\left(k+1\right)}{2^{2k}k!}\sum_{k_{k-1}=0}^{n}\ldots\sum_{k_2=0}^{k_3}\sum_{k_1=0}^{k_2}\left(-1\right)^{k_n}\left(-2\right)^{k_2-k_1}\ldots\left(-k\right)^{n-k_{k-1}}, \qquad n\in \mathbb{N}_0.$$
\end{lemma}

Furthermore, the next lemma provides the asymptotic form for the square of the radius of convexity of the normalized Bessel functions $g_\nu(z)$. The proof of this lemma uses some results of \cite{BPS14}.

\begin{lemma}\label{Lemma_g_asy}
For $\nu>-1$ the radius of convexity $r^c(g_\nu)$ of the function
$$z\mapsto g_\nu(z) =2^\nu\Gamma(\nu+1) z^{1-\nu}J_\nu(z)$$
has the asymptotic behavior
$$\left(r^c\left(g_\nu\right)\right)^2=\nu\left(c+\mathcal{O}\left(\frac{1}{\nu}\right)\right),$$
as $\nu\to \infty$, where $c$ is some positive constant.
\end{lemma}

The next lemma provides the asymptotic form for the radius of convexity of the normalized Bessel functions $h_\nu(z)$. In the proof we use analogous results of \cite[Theorem 1]{BPS14} for Dini function $e_\nu(z)$.

\begin{lemma}\label{Lemma4}
For $\nu>-1$ the radius of convexity $r^c(h_\nu)$ of the function
$$z\mapsto h_\nu(z) =2^\nu\Gamma(\nu+1) z^{1-\frac{\nu}{2}}J_\nu(\sqrt{z})$$
has the asymptotic behavior
$$r^c\left(h_\nu\right)=\nu\left(d+\mathcal{O}\left(\frac{1}{\nu}\right)\right),$$
as $\nu\to \infty$, where $d$ is some positive constant.
\end{lemma}

Before we state the main theorems of this paper, let us define the ordinary potential polynomials. Let $f(z)=1+\sum_{n=1}^{\infty}a_nz^n$ be a formal power series. Corresponding to $f(z)$, for any complex number $\alpha$, the ordinary potential polynomial $A_{\alpha,n}(a_1,a_2,...,a_n)$ is defined by the generating function
$$\left(f(z)\right)^\alpha=\left(1+\sum_{n=0}^{\infty}a_nz^n\right)^\alpha=\sum_{n=0}^{\infty}A_{\alpha,n}(a_1,\dotsc,a_n)z^n.$$
Thus, specifically $A_{\alpha,0}=1,A_{\alpha,1}=\alpha a_1$ and $A_{\alpha,2}=\alpha a_2+\binom{\alpha}{2} a_1^2$. One can refer to \cite{Ne13} for additional details about the ordinary potential polynomials. The next two theorems provide asymptotic expansions for the radius of convexity for two types of normalized Bessel functions of the first kind. The idea of the proofs of the next theorems is inspired from \cite{BN21}.

\begin{theorem}\label{Theorem1}
Let $\eta_n^{(k)}$ denote the coefficients of the expansion in \eqref{etaseries}. Then, the square of the radius of convexity $r^c\left(g_\nu\right)$ has the asymptotic expansion
\begin{equation}\label{asy_relation_g}
\left(r^c\left(g_\nu\right)\right)^2 \sim  \nu\left(c+\sum_{n=1}^{\infty}\frac{\epsilon_n}{\nu^n}\right)
\end{equation}	
as $\nu\to \infty$, where the coefficients $\epsilon_n$ can be determined by the recurrence relation
\begin{equation}\label{epsilon_rln}
\frac{3}{4}\left((-1)^{n+1}c+\sum_{m=0}^{n}(-1)^{n-m}\epsilon_{m+1}\right)+\sum_{k=0}^{n+1}\left(\sum_{m=1}^{\infty}A_{m+1,k}\left(\epsilon_1,\ldots,\epsilon_k\right)\eta_{n-k+1}^{(m+1)}\right)=0
\end{equation}
and $c$ satisfies
\begin{equation}\label{c_value}
c=\frac{2}{3}-\frac{4}{3}\sum_{m=1}^{\infty}A_{m+1,0}\eta_0^{(m+1)}.
\end{equation}
In particular, for $n=0$ in equation \eqref{epsilon_rln} we arrive at
\begin{equation}\label{eps_1_value}
\epsilon_1=\frac{\frac{3c}{4}-\sum_{m=1}^{\infty}c^{m+1}\eta_1^{(m+1)}}{\frac{3}{4}+\sum_{m=1}^{\infty}(m+1)c^m\eta_0^{(m+1)}}.
\end{equation}
By using Mathematica, the coefficients, accurate up to $10^{-5}$, in the above asymptotic expansion are
\begin{equation}\label{conv_num_val}
	\left(r^{c}\left(g_\nu\right)\right)^2 \sim  \nu\left(0.535898+\frac{0.335953}{\nu}+\ldots\right).
\end{equation}
\end{theorem}

\begin{figure}[t]
	\centering
	\includegraphics[width=0.6\textwidth]{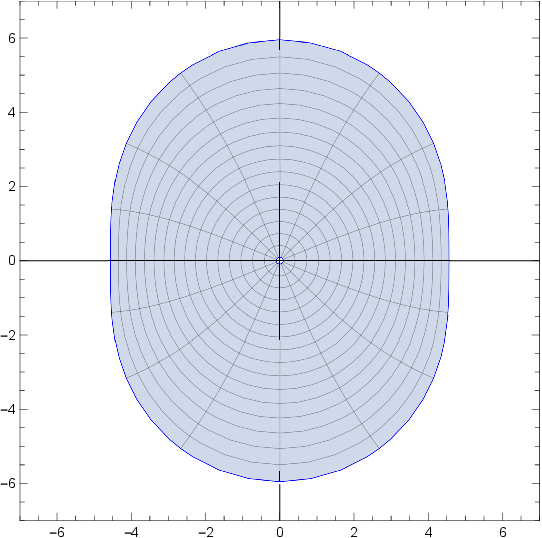}
	\caption{The image of the open disk $\mathbb{D}_{r}$ under the Bessel function $z\mapsto g_\nu(z),$ where $r\sim 5.208\ldots$ is the approximative value of the radius of convexity of $g_\nu(z)$ considering the first two terms of \eqref{conv_num_val} for $\nu=50$.}
	\label{Fig1}
\end{figure}

For large $\nu=50$, by using the first two terms in \eqref{conv_num_val}, we calculated the approximative value of the radius of convexity of $g_\nu(z)$ and ploted in Figure \ref{Fig1}.

\begin{remark}\label{remark1_c}
{\em We would like to mention that we can write \eqref{epsilon_rln} explicitly to determine $\epsilon_n$ for $n>1$. From the proof of Theorem \ref{Theorem1} the expression $A_{k,n}(\epsilon_1,\ldots,\epsilon_n)$ is a potential polynomial given by the generating function
\begin{align*}
\left(c+\sum_{n=1}^{\infty}\frac{\epsilon_n}{\nu^n}\right)^{k}=\sum_{n=0}^{\infty}\frac{A_{k,n}\left(\epsilon_1,\ldots,\epsilon_n\right)}{\nu^n}.
\end{align*}
The above equation can be written as
\begin{equation}\label{r1_A_exp}
c^k\left(1+\sum_{n=1}^{\infty}\frac{\epsilon_n}{c}\frac{1}{\nu^n}\right)^k=\sum_{n=0}^{\infty}A_{k,n}(\epsilon_1,\ldots,\epsilon_n)\frac{1}{\nu^n}.
\end{equation}
Moreover, in view of \cite[Appendix]{Ne13} we have
\begin{equation}\label{r1_tilde_A_exp}
c^k\left(1+\sum_{n=1}^{\infty}\frac{\epsilon_n}{c}\frac{1}{\nu^n}\right)^k=c^k\sum_{n=0}^{\infty}\tilde{A}_{k,n}\left(\frac{\epsilon_1}{c},\ldots,\frac{\epsilon_n}{c}\right)\frac{1}{\nu^n},
\end{equation}
where
$$\tilde{A}_{k,n}\left(\frac{\epsilon_1}{c},\ldots,\frac{\epsilon_n}{c}\right)=\sum_{}^{}\binom{k}{l}\frac{l!}{l_1! l_2!\ldots l_n!}\epsilon_1^{l_1}\epsilon_2^{l_2}\ldots\epsilon_n^{l_n}\frac{1}{c^l}$$
and the sum extends over all sequences $l_1,l_2,\ldots,l_n$ of non-negative integers such that $l_1+2l_2+\ldots+nl_n=n$ and $l_1+l_2+\ldots+l_n=l$.  It is worth to note that for any value of $l$,
$l_n=0\text{ or }1$. Particularly, for $n\in \mathbb{N}$, $l_n = 1$ if and only if $l = 1$. In view of this observation,  the above equation  can be written as 	
\begin{equation}\label{r1_tA_epsilon}
\tilde{A}_{k,n}\left(\frac{\epsilon_1}{c},\ldots,\frac{\epsilon_n}{c}\right)=k\frac{\epsilon_n}{c}+\sum_{}^{}\binom{k}{l}\frac{l!}{l_1! l_2!\ldots l_n!}\epsilon_1^{l_1}\epsilon_2^{l_2}\ldots\epsilon_{n-1}^{l_{n-1}}\frac{1}{c^l},
\end{equation}
where each $l_i$ satisfies the above condition along with $l>1$. From equations \eqref{r1_A_exp}, \eqref{r1_tilde_A_exp} and \eqref{r1_tA_epsilon} we obtain that
\begin{equation}\label{r1_exp_A}
A_{k,n}(\epsilon_1,\ldots,\epsilon_n)=c^k\left(k\frac{\epsilon_n}{c}+\sum_{}^{}\binom{k}{l}\frac{l!}{l_1!\ldots l_n!}\epsilon_1^{l_1}\epsilon_2^{l_2}\ldots\epsilon_{n-1}^{l_{n-1}}\frac{1}{c^l}\right).
\end{equation}
Moreover, we can rewrite equation \eqref{epsilon_rln} after separating the terms associated with $\epsilon_{n+1}$ as
\begin{align*}
&\frac{3}{4}\left((-1)^{n+1}c+\sum_{m=0}^{n-1}(-1)^{n-m}\epsilon_{m+1}\right)+\frac{3}{4}\epsilon_{n+1}+\sum_{k=0}^{n}\left(\sum_{m=1}^{\infty}A_{m+1,k}\left(\epsilon_1,\ldots,\epsilon_k\right)\eta_{n-k+1}^{(m+1)}\right)\\
&\qquad+\sum_{m=1}^{\infty}A_{m+1,n+1}\left(\epsilon_1,\ldots,\epsilon_{n+1}\right)\eta_0^{(m+1)}=0
\end{align*}	
Now, putting the value of $A_{m+1,n+1}(\epsilon_1,\ldots,\epsilon_{n+1})$ by using \eqref{r1_exp_A} in the above equation and solving for $\epsilon_{n+1}$, we obtain that	
\begin{align*}
\epsilon_{n+1}&\left(\frac{3}{4}+ \sum_{m=1}^{\infty}\eta_{0}^{(m+1)}(m+1)c^{m}\right)\\&=-\sum_{m=1}^{\infty}c^{m+1}\eta_{0}^{(m+1)}\sum_{}^{}\binom{m+1}{l}\frac{l!}{l_1!\ldots l_{n+1}!}\epsilon_1^{l_1}\epsilon_2^{l_2}\ldots\epsilon_{n}^{l_{n}}\frac{1}{c^l}\\
&-\frac{3}{4}\left((-1)^{n+1}c+\sum_{m=0}^{n-1}(-1)^{n-m}\epsilon_{m+1}\right)-\sum_{k=0}^{n}\left(\sum_{m=1}^{\infty}A_{m+1,k}\left(\epsilon_1,\ldots,\epsilon_k\right)\eta_{n-k+1}^{(m+1)}\right),
\end{align*}
where $l>1$, $l_1+2l_2+\ldots+nl_n=n$ and $l_1+l_2+\ldots+l_n=l$. Notice that the condition $l>1$ results in $l_{n+1}=0$.}
\end{remark}

\begin{remark}\label{polyn_c}
{\em We also note that from the proof of Theorem \ref{Theorem1} we have that $A_{m+1,0}=c^{m+1}$. Hence \eqref{c_value} can be written as
\begin{equation}\label{c_series}
c=\frac{2}{3}-\frac{4}{3}\sum_{m=1}^{\infty}c^{m+1}\eta_0^{(m+1)}.
\end{equation}
By using the fact that the denominator of the terms in the series $\eta_0^{(k)}$ increases rapidly, the infinite series \eqref{c_series} can be reduced into a polynomial of a finite yet large degree in order to numerically approximate the value of $c$, which satisfies the series \eqref{c_series}.
For example the $20^{th}$ term of the series \eqref{conv_num_val} has value less than $10^{-7}$. We find out the value of $c$ by solving the polynomial of degree $20$ with the variable $c$, which is given in \eqref{conv_num_val}.}
\end{remark}

\begin{remark}
{\em By using the Euler-Rayleigh inequalities, Akta\c{s} et al. \cite[Theorem 6]{ABO18} derived the bounds for the radius of convexity of normalized Bessel function $g_\nu(z)$ as follows
\begin{equation}\label{ineq_g}
4\sqrt{\frac{(\nu+1)^2\left(\nu+2\right)}{56\nu+137}}<\left(r^c(g_\nu)\right)^2<\frac{2\left(56\nu+137\right)\left(\nu+1\right)\left(\nu+3\right)}{208\nu^2+1172\nu+1693}.
\end{equation}
The right-hand side of \eqref{ineq_g} can be expressed as
\begin{align*} \frac{2\left(56\nu+137\right)\left(\nu+1\right)\left(\nu+3\right)}{208\nu^2+1172\nu+1693}&=\frac{7\nu\left(1+\frac{137}{56\nu}\right)\left(1+\frac{1}{\nu}\right)\left(1+\frac{3}{\nu}\right)}{13\left(1+\frac{1172}{208\nu}+\frac{1693}{208\nu^2}\right)}=\nu\left(\frac{7}{13}+\frac{591}{1352}\frac{1}{\nu}+\mathcal{O}\left(\frac{1}{\nu^2}\right)\right).
\end{align*}
Similarly, we can express the left-hand side of the inequality \eqref{ineq_g} as
\begin{align*}
4\sqrt{\frac{(\nu+1)^2\left(\nu+2\right)}{56\nu+137}}=\nu\left(\sqrt{\frac{2}{7}}+\frac{87}{56\sqrt{14}}\frac{1}{\nu}+\mathcal{O}\left(\frac{1}{\nu^2}\right)\right).
\end{align*}
In view of \eqref{asy_relation_g} and the above discussion, we can say that $c$ lies in the interval $\left(\sqrt{{2}/{7}},{7}/{13}\right)$.
This interval could be further narrowed by taking higher order  Euler-Rayleigh inequalities (see the proof of \cite[Theorem 6]{ABO18}). The value of $c$, which serves as the root of the polynomial mentioned in Remark \ref{polyn_c} and lies within the precise interval, can be regarded as an approximate value for the constant $c$ that fits into the asymptotic expansion of the square of $r^c(g_\nu)$. It is worth also to note that the accuracy of the value of $c$ can be improved by narrowing the interval and using a higher degree polynomial.}
\end{remark}

\begin{remark}
{\em It is worth mentioning that besides using the recurrence relation \eqref{epsilon_rln} to calculate $\epsilon_1,\epsilon_2,\ldots$, for large $\nu$, we can find their bounds by using the  Euler-Rayleigh inequalities.
In view of relation \eqref{asy_relation_g} and inequality \eqref{ineq_g}, which bounds the radius of convexity $r^c\left(g_\nu\right)$, for large $\nu$ we obtain a bound for $\epsilon_1$ as
$$|\epsilon_1|<r,$$
where
\begin{equation}\label{def_r}
	r=\frac{591}{1352}=\max\left(\left|\frac{591}{1352}\right|,\left|\frac{87}{56\sqrt{14}}\right|\right).
\end{equation}
In a similar way we can find bounds for $\epsilon_2,\epsilon_3,\ldots$, for large $\nu$.
For example it is clear from \eqref{def_r} that $r=0.4371302$ and
$$|\epsilon_1|<0.4371302.$$
This bound is satisfied by the value we calculated by using the recurrence relation \eqref{epsilon_rln} which is given $0.335953$ (see \eqref{conv_num_val}).
From the proof of the Lemma \ref{Lemma_g_asy}, we noted that the bounds of the radius of convexity converge, so we can find tighter bounds for $\epsilon_1$ by using other bounds for the radius of convexity (see also the proof of \cite[Theorem 6]{ABO18}). Note that we used modulus in \eqref{def_r} to emphasize that the coefficients could be negative for other Euler-Rayleigh inequalities.}
\end{remark}

Before stating the next theorem, let us note that the approach used in both of the previous remarks is also applicable for the following theorem to approximate the constant $d$ and finding bounds for the coefficients $\epsilon_1,\epsilon_2,\ldots$ in the asymptotic expansion of $r^c\left(h_\nu\right)$.

\begin{theorem}\label{Theorem2}
Let $\theta_n^{(k)}$ denote the coefficients of the expansion in \eqref{thetaseries}. Then, the radius of convexity $r^c\left(h_\nu\right)$, of the function $h_\nu(z)$ defined in \eqref{unif_h_nu_def}, has the asymptotic expansion
\begin{align*}
r^c\left(h_\nu\right) \sim  \nu\left(d+\sum_{n=1}^{\infty}\frac{\epsilon_n}{\nu^n}\right)
\end{align*}	
as $\nu\to \infty$, where the coefficients $\epsilon_n$ can be determined by the recurrence relation
\begin{equation}\label{epsilon_rln_h}
\frac{3}{4}\left((-1)^{n+1}d+\sum_{m=0}^{n}(-1)^{n-m}\epsilon_{m+1}\right)+\sum_{k=0}^{n+1}\left(\sum_{m=1}^{\infty}A_{m+1,k}\left(\epsilon_1,\ldots,\epsilon_k\right)\theta_{n-k+1}^{(m+1)}\right)=0
\end{equation}
and $d$ satisfies
$$d=\frac{4}{3}\left(1-\sum_{m=1}^{\infty}A_{m+1,0}\theta_0^{(m+1)}\right).$$
In particular, for $n=0$ in equation \eqref{epsilon_rln_h} we obtain that
\begin{align*}
\epsilon_1=\frac{\frac{3d}{4}-\sum_{m=1}^{\infty}d^{m+1}\theta_1^{(m+1)}}{\frac{3}{4}+\sum_{m=1}^{\infty}(m+1)d^m\theta_0^{(m+1)}}.
\end{align*}
By using Mathematica, the coefficients, accurate up to $10^{-5}$, in the above asymptotic expansion are
\begin{equation}\label{h_conv_num_val}
	r^{c}\left(h_\nu\right) \sim  \nu\left(1.17157+\frac{0.858757}{\nu}+\ldots\right).
\end{equation}
\end{theorem}

\begin{figure}[t]
	\centering
	\includegraphics[width=0.6\textwidth]{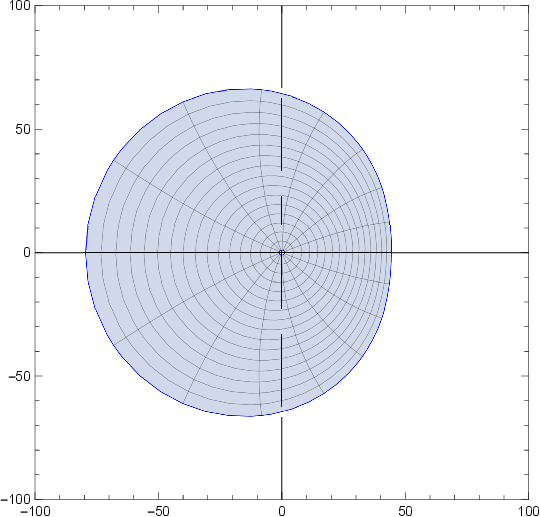}
	\caption{The image of the open disk $\mathbb{D}_{r}$ under the Bessel function $z\mapsto h_\nu(z),$ where $r\sim 59.437\ldots$ is the approximative value of the radius of convexity of $h_\nu(z)$ considering the first two terms of \eqref{h_conv_num_val} for $\nu=50$.}
	\label{Fig2}
\end{figure}

For large $\nu=50$, by using the first two terms in \eqref{h_conv_num_val}, we calculated the approximative value of the radius of convexity of $h_\nu(z)$ and we have shown this in Figure \ref{Fig2}.

\subsection{Bounds and asymptotic expansions for the radii of uniform convexity of normalized Bessel functions}
Now, we focus on similar results as before related to uniform convexity. The next lemmas provide bounds for the radii of uniform convexity of normalized Bessel functions of the first kind. We use the so-called Euler-Rayleigh inequality technique \cite{IM95} and properties of the Laguerre-P\'olya class of entire functions to deduce these bounds.

\begin{lemma}\label{unif_bounds_g}
	Let $\nu>-1$. Then the radius of uniform convexity $r^{uc}(g_\nu)$ of the function
	$$z\mapsto g_\nu(z)=2^\nu\Gamma(\nu+1)z^{1-\nu}J_\nu(z)$$
	is the smallest positive root of the equation
	\begin{align*}
		g_\nu^{\prime}(z)+2 z g_\nu^{\prime\prime}(z)=0
	\end{align*}
	and satisfies the following inequality
	\begin{align*}
		\omega_k^{-\frac{1}{k}}<(r^{uc}(g_\nu))^2<\frac{\omega_k}{\omega_{k+1}},
	\end{align*}
	where
	$\omega_k$ is given by \eqref{unif_bound_exprsn}.
	In particular for $k=1$ we obtain that
	\begin{equation}\label{unif_first_inq}
		2\sqrt{\frac{(\nu+1)}{15}}<r^{uc}(g_\nu)<2\sqrt{\frac{\nu(\nu+1)}{3(4\nu-1)}}.
	\end{equation}
\end{lemma}

\begin{lemma}\label{unif_bounds_h}
	Let $\nu>-1$. Then the radius of uniform convexity $r^{uc}(h_\nu)$ of the function
	$$z\mapsto h_\nu(z)=2^\nu\Gamma(\nu+1) z^{1-\frac{\nu}{2}}J_\nu(\sqrt{z})$$
	is the smallest positive root of the equation
	\begin{align*}
		h_\nu^{\prime}(z)+2 z h_\nu^{\prime\prime}(z)=0
	\end{align*}
	and satisfies the following inequality
	\begin{align*}
		\sigma_k^{-\frac{1}{k}}<r^{uc}(h_\nu)<\frac{\sigma_k}{\sigma_{k+1}},
	\end{align*}
	where
	$\sigma_k$ is given by \eqref{unif_bound_exprsn_h}.
\end{lemma}

Furthermore, the next lemma provides the asymptotic form for the square of the radius of uniform convexity of the normalized Bessel functions $g_\nu(z)$.

\begin{lemma}\label{unif_Lemma_g_asy}
	The radius of uniform convexity $r^{uc}(g_\nu)$ of the function
	$$z\to g_\nu(z) =2^\nu\Gamma(\nu+1) z^{1-\nu}J_\nu(z)$$
	admits the asymptotic behavior
	$$\left(r^{uc}\left(g_\nu\right)\right)^2=\nu\left(c+\mathcal{O}\left(\frac{1}{\nu}\right)\right),$$
	for $\nu\to \infty$, where $c$ is some positive constant.
\end{lemma}

Similarly to the previous lemma, the next lemma provides the asymptotic form for the radius of uniform convexity of the normalized Bessel functions $h_\nu(z)$.
\begin{lemma}\label{unif_Lemma_h_asy}
	The radius of uniform convexity $r^{uc}(h_\nu)$ of the function
	$$z\to h_\nu(z) =2^\nu\Gamma(\nu+1) z^{1-\frac{\nu}{2}}J_\nu(\sqrt{z})$$
	admits the asymptotic behavior
	$$r^{uc}\left(h_\nu\right)=\nu\left(d+\mathcal{O}\left(\frac{1}{\nu}\right)\right),$$
	for $\nu\to \infty$, where $d$ is some positive constant.
\end{lemma}

The next two theorems are the main results of this subsection. The idea of the proofs of the next theorems is also inspired from \cite{BN21}.

\begin{theorem}\label{unif_Theorem1}
	Let $\eta_n^{(k)}$ denote the coefficients of the expansion in \eqref{etaseries}. Then, the square of the radius of uniform convexity $r^{uc}\left(g_\nu\right)$ has the asymptotic expansion
	\begin{equation}\label{unif_asy_relation_g}
		\left(r^{uc}\left(g_\nu\right)\right)^2 \sim  \nu\left(\tilde{c}+\sum_{n=1}^{\infty}\frac{\varepsilon_n}{\nu^n}\right)
	\end{equation}	
	as $\nu\to \infty$, where the coefficients $\varepsilon_n$ can be determined by the recurrence relation
	\begin{equation}\label{unif_epsilon_rln}
		\frac{3}{4}\left((-1)^{n+1}\tilde{c}+\sum_{m=0}^{n}(-1)^{n-m}\varepsilon_{m+1}\right)+\sum_{k=0}^{n+1}\left(\sum_{m=1}^{\infty}A_{m+1,k}\left(\varepsilon_1,\ldots,\varepsilon_k\right)\eta_{n-k+1}^{(m+1)}\right)=0
	\end{equation}
	and $\tilde{c}$ satisfies
	\begin{equation}\label{unif_c_value}
		\tilde{c}=\frac{1}{3}-\frac{4}{3}\sum_{m=1}^{\infty}A_{m+1,0}\eta_0^{(m+1)}.
	\end{equation}
	In particular, for $n=0$ in equation \eqref{unif_epsilon_rln} we arrive at
	\begin{align*}
		\varepsilon_1=\frac{\frac{3 \tilde{c}}{4}-\sum_{m=1}^{\infty}\tilde{c}^{m+1}\eta_1^{(m+1)}}{\frac{3}{4}+\sum_{m=1}^{\infty}(m+1)\tilde{c}^m\eta_0^{(m+1)}}.
	\end{align*}
	By using Mathematica, the coefficients, accurate up to $10^{-5}$, in the above asymptotic expansion are
	\begin{equation}\label{unif_num_val}
		\left(r^{uc}\left(g_\nu\right)\right)^2 \sim  \nu\left(0.298438+\frac{0.218612}{\nu}+\ldots\right).
	\end{equation}
\end{theorem}

\begin{figure}[t]
	\centering
	\includegraphics[width=0.6\textwidth]{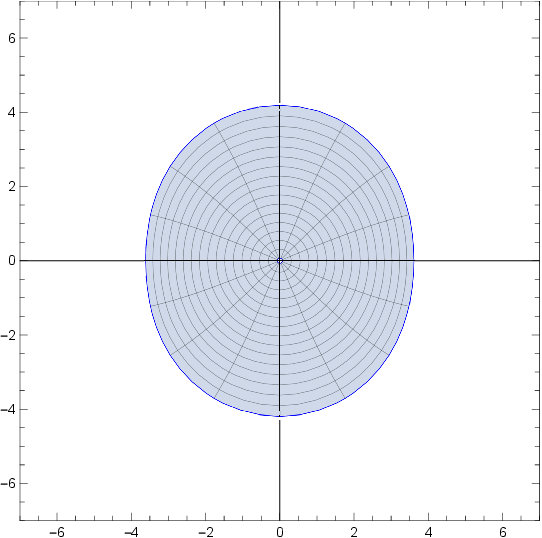}
	\caption{The image of the open disk $\mathbb{D}_{r}$ under the Bessel function $z\mapsto g_\nu(z),$ where $r\sim 3.891\ldots$ is the approximative value of the radius of uniform convexity of $g_\nu(z)$ considering the first two terms of \eqref{unif_num_val} for $\nu=50$.}
	\label{Fig3}
\end{figure}

For large $\nu=50$, by using the first two terms in \eqref{unif_num_val}, we calculated the approximative value of the radius of uniform convexity of $g_\nu(z)$ and ploted in Figure \ref{Fig3}.

\begin{remark}
	{\em By using a similar argument as in Remark \ref{remark1_c}, we can write \eqref{unif_epsilon_rln} explicitly to determine $\varepsilon_n$ for $n>1$ as
		\begin{align*}
			\varepsilon_{n+1}&\left(\frac{3}{4}+ \sum_{m=1}^{\infty}\eta_{0}^{(m+1)}(m+1)\tilde{c}^{m}\right)\\&=-\sum_{m=1}^{\infty}\tilde{c}^{m+1}\eta_{0}^{(m+1)}\sum_{}^{}\binom{m+1}{l}\frac{l!}{l_1!\ldots l_{n+1}!}\varepsilon_1^{l_1}\varepsilon_2^{l_2}\ldots\varepsilon_{n}^{l_{n}}\frac{1}{\tilde{c}^l}\\
			&-\frac{3}{4}\left((-1)^{n+1}\tilde{c}+\sum_{m=0}^{n-1}(-1)^{n-m}\varepsilon_{m+1}\right)-\sum_{k=0}^{n}\left(\sum_{m=1}^{\infty}A_{m+1,k}\left(\varepsilon_1,\ldots,\varepsilon_k\right)\eta_{n-k+1}^{(m+1)}\right)
		\end{align*}
		Here, $l>1$, $l_1+2l_2+\ldots+nl_n=n$ and $l_1+l_2+\ldots+l_n=l$. Notice that the condition $l>1$ results in $l_{n+1}=0$.}
\end{remark}

\begin{remark}\label{unif_polyn_c}
	{\em Note that from the proof of Theorem \ref{unif_Theorem1}, $A_{m+1,0}=\tilde{c}^{m+1}$. Hence \eqref{unif_c_value} can be written as
		\begin{equation}\label{unif_c_series}
			\tilde{c}=\frac{1}{3}-\frac{4}{3}\sum_{m=1}^{\infty}\tilde{c}^{m+1}\eta_0^{(m+1)}.
		\end{equation}
		By using Lemma \ref{lemma_eta}, we observed that $\eta_0^{m+1}$ decreases rapidly for large $m$,
		%	By using  Lemma 1, we observed that $\eta_0^(m+1}$ decreases rapidly for large $m$,	
	the infinite series \eqref{unif_c_series} can be reduced into a polynomial of a finite yet large degree in order to numerically approximate the value of $\tilde{c}$, which satisfies the series \eqref{unif_c_series}. For example, the $20^{th}$ term of the series \eqref{unif_c_series} has value less than $10^{-7}$. We find out the value of $\tilde{c}$ by solving the polynomial of degree $20$ with the variable $\tilde{c}$, which is given in \eqref{unif_num_val}.}
	\end{remark}
	
	\begin{remark}
{\em From Lemma \ref{unif_bounds_g} we have that
	\begin{equation}\label{unif_ineq_g}
		\frac{4(\nu+1)}{15}<\left(r^{uc}(g_\nu)\right)^2<\frac{4\nu(\nu+1)}{3(4\nu-1)}.
	\end{equation}
	The right-hand side of \eqref{unif_ineq_g} can be expressed as
	\begin{align*} \frac{4\nu(\nu+1)}{3(4\nu-1)}&=\frac{\nu}{3}\left(1 +\frac{5}{4\nu} + \mathcal{O}\left(\frac{1}{\nu^2}\right)\right).
	\end{align*}
	Similarly, we can express the left-hand side of the inequality \eqref{unif_ineq_g} as
	\begin{align*}
		\frac{4(\nu+1)}{15}=\frac{4\nu}{15}\left(1+\frac{1}{\nu}\right)
	\end{align*}
	In view of \eqref{unif_asy_relation_g} and the above discussion, we can say that $\tilde{c}$ lies in the interval $\left(\frac{4}{15},\frac{1}{3}\right)$.
	This interval could be further narrowed by taking higher order  Euler-Rayleigh inequalities (see the proof of Lemma \ref{unif_bounds_g}). The value of $\tilde{c}$, which serves as the root of the polynomial mentioned in Remark \ref{unif_polyn_c} and lies within the precise interval, can be regarded as an approximate value for the constant $\tilde{c}$ that fits into the asymptotic expansion of the square of $r^{uc}(g_\nu)$.
	It is noteworthy that the accuracy of the value of $\tilde{c}$ can be improved by narrowing the interval and using a higher degree polynomial.}
	\end{remark}
	
	\begin{remark}
{\em It is worth mentioning that besides using recurrence relation \eqref{unif_epsilon_rln} to calculate $\varepsilon_1,\varepsilon_2,\ldots$, for large $\nu$, we can also find their bounds by using the  Euler–Rayleigh inequalities.
	In view of relation \eqref{unif_asy_relation_g} and inequality \eqref{unif_ineq_g}, which bounds the radius of uniform convexity $r^{uc}\left(g_\nu\right)$, for large $\nu$ we obtain a bound for $\varepsilon_1$ as
	$$|\varepsilon_1|<r,$$
	where
	\begin{equation}\label{unif_def_r}
		r=\frac{5}{12}=\max\left(\left|\frac{4}{15}\right|,\left|\frac{5}{12}\right|\right).
	\end{equation}
	In a similar way we can find bounds for $\varepsilon_2,\varepsilon_3,\ldots$, for large $\nu$. For example it is clear from \eqref{unif_def_r} that $r=0.41666$ and
	$$|\varepsilon_1|<0.41666.$$
	This bound is satisfied by the value we calculated by using the recurrence relation \eqref{unif_epsilon_rln} which is given $0.218612$ (see \eqref{unif_num_val}).
	From the proof of the Lemma \ref{unif_Lemma_g_asy}, we noted that the bounds of the radius of uniform convexity converge, so we can find tighter bounds for $\varepsilon_1$ by using other bounds for the radius of uniform convexity (see also the proof of Lemma \ref{unif_bounds_g}). Note that we used modulus in \eqref{unif_def_r} to emphasize that the coefficients could be negative for other Euler-Rayleigh inequalities.}
	\end{remark}
	
	Before stating the next theorem, let us note that the approach used in both of the previous remarks is also applicable for the following theorem to approximate the constant $\tilde{d}$ and finding bounds for the coefficients $\varepsilon_1,\varepsilon_2,\ldots$ in the asymptotic expansion of $r^{uc}\left(h_\nu\right)$.
	
\begin{theorem}\label{unif_Theorem2}
Let $\theta_n^{(k)}$ denote the coefficients of the expansion in \eqref{thetaseries}. Then, the radius of uniform convexity $r^{uc}\left(h_\nu\right)$ has the asymptotic expansion
\begin{align*}
	r^{uc}\left(h_\nu\right) \sim  \nu\left(\tilde{d}+\sum_{n=1}^{\infty}\frac{\varepsilon_n}{\nu^n}\right)
\end{align*}	
as $\nu\to \infty$, where the coefficients $\varepsilon_n$ can be determined by the recurrence relation
\begin{equation}\label{unif_epsilon_rln_h}
	\frac{3}{4}\left((-1)^{n+1}\tilde{d}+\sum_{m=0}^{n}(-1)^{n-m}\varepsilon_{m+1}\right)+\sum_{k=0}^{n+1}\left(\sum_{m=1}^{\infty}A_{m+1,k}\left(\varepsilon_1,\ldots,\varepsilon_k\right)\theta_{n-k+1}^{(m+1)}\right)=0
\end{equation}
and $\tilde{d}$ satisfies
$$
\tilde{d}=\frac{2}{3}-\frac{4}{3}\sum_{m=1}^{\infty}A_{m+1,0}\theta_0^{(m+1)}.
$$
In particular, for $n=0$ in equation \eqref{unif_epsilon_rln_h} we obtain that
\begin{align*}
	\varepsilon_1=\frac{\frac{3 \tilde{d}}{4}-\sum_{m=1}^{\infty}\tilde{d}^{m+1}\theta_1^{(m+1)}}{\frac{3}{4}+\sum_{m=1}^{\infty}(m+1)\tilde{d}^m\theta_0^{(m+1)}}.
\end{align*}
By using Mathematica, the coefficients, accurate up to $10^{-5}$, in the above asymptotic expansion are
\begin{equation}\label{h_unif_conv_num_val}
	r^{c}\left(h_\nu\right) \sim  \nu\left(0.627719+\frac{0.478612}{\nu}+\ldots\right).
\end{equation}
\end{theorem}

\begin{figure}[t]
	\centering
	\includegraphics[width=0.6\textwidth]{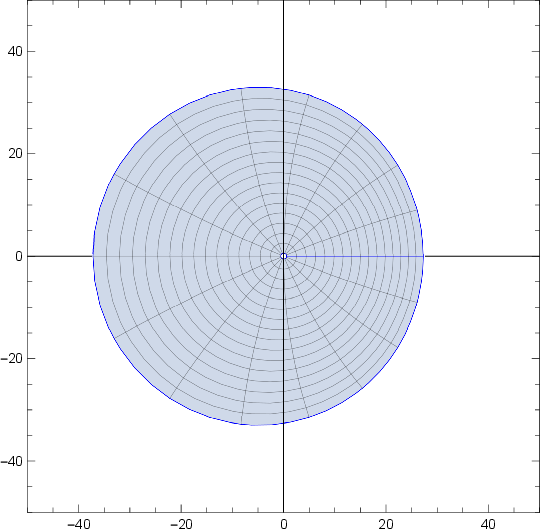}
	\caption{The image of the open disk $\mathbb{D}_{r}$ under the Bessel function $z\mapsto h_\nu(z),$ where $r\sim 31.86\ldots$ is the approximative value of the radius of uniform convexity of $h_\nu(z)$ considering the first two terms of \eqref{h_unif_conv_num_val} for $\nu=50$.}
	\label{Fig4}.
\end{figure}

For large $\nu=50$, by using the first two terms in \eqref{h_unif_conv_num_val}, we calculated the approximative value of the radius of uniform convexity of $h_\nu(z)$ and ploted in Figure \ref{Fig4}.

\section{Proofs of the preliminary and main results}\label{proofs}
\setcounter{equation}{0}

\begin{proof}[\bf Proof of Lemma \ref{lemma_eta}]
Due to \cite[Theorem 1]{BPS14} for $\nu>-1$ the Dini function $d_\nu$, as defined in \eqref{dini_d}, admits the Weistrassian factorization as
\begin{align*}
d_\nu(z)=\frac{z^\nu}{2^\nu\Gamma\left(\nu+1\right)}\prod_{n=1}^{\infty}\left(1-\frac{z^2}{\alpha_{\nu,n}^2}\right),
\end{align*}
where the infinite product is uniformly convergent on each compact subset of the complex plane. By replacing $z$ by $\sqrt{z}$ we write the above expression as
$$r_\nu(z)=2^\nu\Gamma(\nu+1)z^{-\frac{\nu}{2}}d_\nu(\sqrt{z})=\prod_{n=1}^{\infty}\left(1-\frac{z}{\alpha_{\nu,n}^2}\right).$$
On the other hand, by using the infinite sum representation of the Bessel function of the first kind \eqref{Bessel_series} and equation \eqref{dini_d}, we obtain
$$r_\nu(z)=1+\sum_{n=1}^{\infty}\frac{(-1)^n\left(2n+1\right)}{2^{2n}n!(\nu+1)\ldots(\nu+n)}z^n=1+\sum_{n=1}^{\infty}a_nz^n.$$
By using the definition of $\eta_k$ and equations \cite[eqns. (3.4) and (3.7)]{IM95} we have the next recurrence relation for $\eta_k(\nu)$
\begin{equation}\label{eta_1}
\eta_1(\nu)=\sum_{m=1}^{\infty}\frac{1}{\alpha_{\nu,n}^{2}}=-a_1
\end{equation}
and
\begin{equation}\label{eta_k}
\eta_k(\nu)=-na_k-\sum_{i=1}^{k-1}a_i\eta_{k-i}(\nu).
\end{equation}

Now, we prove by induction on $k$, that for any positive integer $k$ and positive real $\nu>k$, the Rayleigh sum $\eta_k\left(\nu\right)$, can be written in the form \eqref{etaseries}. First notice that for $k\in\mathbb{N}$
\begin{align*}
a_k&=\frac{\left(-1\right)^k\left(2k+1\right)}{2^{2k}k!\left(\nu+1\right)\left(\nu+2\right)\ldots\left(\nu+k\right)}\\
&=\frac{\left(-1\right)^k\left(2k+1\right)}{\nu^k2^{2k}k!\left(1+\frac{1}{\nu}\right)\left(1+\frac{2}{\nu}\right)\ldots\left(1+\frac{k}{\nu}\right)}\\ &=\frac{\left(-1\right)^k\left(2k+1\right)}{\nu^k2^{2k}k!}\left(\sum_{{k_1}=0}^{\infty}\frac{\left(-1\right)^{k_1}}{\nu^{k_1}}\right)\left(\sum_{k_2=0}^{\infty}\frac{\left(-2\right)^{k_2}}{\nu^{k_2}}\right)\ldots\left(\sum_{k_k=0}^{\infty}\frac{\left(-k\right)^{k_k}}{\nu^{k_k}}\right)\\
&=\frac{\left(-1\right)^k\left(2k+1\right)}{\nu^k2^{2k}k!}\sum_{k_k=0}^{\infty}\sum_{k_{k-1}=0}^{k_k}\ldots\sum_{k_2=0}^{k_3}
\sum_{k_1=0}^{k_2}\frac{\left(-1\right)^{k_1}\left(-2\right)^{k_2-k_1}\ldots\left(-k\right)^{k_k-k_{k-1}}}{\nu^{k_k}}
\end{align*}
or equivalently
\begin{equation}\label{a_k_expansion}
a_k=\frac{1}{\nu^k}\sum_{n=0}^{\infty}\frac{a_{n}^{(k)}}{\nu^{n}},
\end{equation}
where            $$a_{n}^{(k)}=\frac{\left(-1\right)^k\left(2k+1\right)}{2^{2k}k!}\sum_{k_{k-1}=0}^{n}\ldots\sum_{k_2=0}^{k_3}\sum_{k_1=0}^{k_2}\left(-1\right)^{k_1}\left(-2\right)^{k_2-k_1}\ldots\left(-k\right)^{n-k_{k-1}}.$$
For $k=1$ from equations \eqref{eta_1} and \eqref{a_k_expansion} we obtain that $$\eta_1(\nu)=-a_1=-\frac{1}{\nu}\sum_{k_1=0}^{\infty}\frac{a^{(1)}_{k_1}}{\nu^{k_n}}=\frac{1}{\nu}\sum_{n=0}^{\infty}\frac{\eta_n^{(1)}}{\nu^n},$$
where $\eta_n^{(1)}=-a_{n}^{(1)}$. Let $N\geq 2$ and suppose that $\eta_{k}(\nu)$ can be expressed in form \eqref{etaseries} for $1\leq k\leq N-1$.
For $k=N$ by using equations \eqref{a_k_expansion}, \eqref{eta_k} and the induction hypothesis we can write
\begin{align*}
\eta_N(\nu)&=-Na_N-\sum_{i=1}^{N-1}a_i\eta_{N-i}(\nu)\\ &=\frac{-N}{\nu^N}\sum_{n=0}^{\infty}\frac{a_n^{(N)}}{\nu^n}-\sum_{i=1}^{N-1}\left(\frac{1}{\nu^i}\sum_{k_i=0}^{\infty}\frac{a_{k_i}^{(i)}}{\nu^{k_i}}\frac{1}{\nu^{N-i}}\sum_{n=0}^{\infty}\frac{\eta_n^{(N-i)}}{\nu^n}\right)\\
&=\frac{1}{\nu^N}\left(-N\sum_{n=0}^{\infty}\frac{a_n^{(N)}}{\nu^n}-\sum_{i=1}^{N-1}\left(\sum_{n=0}^{\infty}\sum_{m=0}^{n}\frac{a_m^{(i)}\eta_{n-m}^{(N-i)}}{\nu^n}\right)\right)\\
&=\frac{1}{\nu^N}\sum_{n=0}^{\infty}\frac{1}{\nu^n}\left(-Na_n^{(N)}-\sum_{m=0}^{n}\sum_{i=1}^{N-1}a_m^{(i)}\eta_{n-m}^{(N-i)}\right)\\
&=\frac{1}{\nu^N}\sum_{n=0}^{\infty}\frac{\eta_n^{(N)}}{\nu^n},
\end{align*}
where $$\eta_n^{(N)}=-Na_n^{(N)}-\sum_{m=0}^{n}\sum_{i=1}^{N-1}a_m^{(i)}\eta_{n-m}^{(N-i)} \quad \text{for all } n\in \mathbb{N}_0.$$
This completes the proof of the Lemma.
\end{proof}

\begin{proof}[\bf Proof of Lemma \ref{lemma_theta}]
For $\nu>-1$ from the proof of \cite[Theorem 1.3]{BS14} we write
$$\frac{h_\nu^{\prime\prime}(z)}{h_\nu^\prime(z)}=-\sum_{n=1}^{\infty}\frac{1}{\beta_{\nu,n}^2-z}.$$
Integrating both sides of the above equation we arrive at
\begin{align*}
\log h_\nu^\prime(z)=\sum_{n=1}^{\infty}\log\left(1-\frac{z}{\beta_{\nu,n}^2}\right)+c_h
\end{align*}
or equivalently at
\begin{equation}\label{h_nu_prod}
h_\nu^\prime(z)=e^{c_h}\prod_{n=1}^{\infty}\left(1-\frac{z}{\beta_{\nu,n}^2}\right),
\end{equation}
where $c_h$ is some constant. On the other hand, differentiating both sides of $h_\nu(z)=2^\nu\Gamma(\nu+1)z^{1-\frac{\nu}{2}}J_\nu(\sqrt{z})$ we have that
\begin{equation}\label{h_nu_prime}
h_\nu^\prime(z)=2^{\nu-1}\Gamma(\nu+1)z^{-\frac{\nu}{2}}e_\nu(\sqrt{z}),
\end{equation}
where $e_\nu(z)$ is given by \eqref{dini_e}. In view of the recurrence relation $zJ_\nu^\prime(z)=-zJ_{\nu+1}(z)+\nu J_\nu(z)$, we obtain that
\begin{equation}\label{e_nu_exp}
e_\nu(z)=2J_\nu(\sqrt{z})-\sqrt{z}J_{\nu+1}(\sqrt{z}).
\end{equation}
Now, by using the infinite sum representation of the Bessel function of the first kind \eqref{Bessel_series} and equations \eqref{h_nu_prime} and \eqref{e_nu_exp}, we obtain $$h_\nu^\prime(z)=2^\nu\Gamma(\nu+1)\sum_{n=0}^{\infty}\frac{(-1)^nz^n}{2^{2n+\nu}n!\Gamma(n+\nu+1)}-2^{\nu-1}\Gamma(\nu+1)\sum_{n=0}^{\infty}\frac{(-1)^nz^{n+1}}{2^{2n+\nu}n!\Gamma(n+\nu+1)},$$ which implies that $h_\nu^\prime(0)=1$. By using \eqref{h_nu_prod}, we arrive at $c_h=0$. Hence the Dini function $e_\nu(z)$, in view of \eqref{h_nu_prime} and \eqref{h_nu_prime}, admits the Weistrassian factorization as
\begin{align*}
e_\nu(z)=\frac{z^\nu}{2^{\nu-1}\Gamma\left(\nu+1\right)}\prod_{n=1}^{\infty}\left(1-\frac{z^2}{\beta_{\nu,n}^2}\right),
\end{align*}
where the infinite product is uniformly convergent on each compact subset of the complex plane. By replacing $z$ by $\sqrt{z}$ we write the above expression as
$$s_\nu(z)=2^{\nu-1}\Gamma(\nu+1)z^{-\frac{\nu}{2}}e_\nu(\sqrt{z})=\prod_{n=1}^{\infty}\left(1-\frac{z}{\beta_{\nu,n}^2}\right).$$
By using \eqref{Bessel_series} together with \eqref{dini_e} we obtain
$$s_\nu(z)=1+\sum_{n=1}^{\infty}\frac{(-1)^n\left(n+1\right)}{2^{2n}n!(\nu+1)\ldots(\nu+n)}z^n=1+\sum_{n=1}^{\infty}b_nz^n,$$
and the remaining steps are quite similar to those of the proof of Lemma \ref{lemma_eta}, so we omit the details.	
\end{proof}

\begin{proof}[\bf Proof of Lemma \ref{Lemma_g_asy}]
	We note that for $\nu>-1,$ $z\in \mathbb{C}$, the function
	$z\mapsto\Delta_\nu(z)=(zg_\nu^\prime(z))^\prime$ has
	the infinite sum  and infinite product representation as follows
	\cite[eq. (2.13) and eq. (2.14)]{ABO18}
	\begin{equation}\label{delta_series}
		\Delta_\nu(z)=1+\sum_{n\geq1}^{}\frac{(-1)^n(2n+1)^2z^{2n}}{2^{2n}n!(\nu+1)_n}
	\end{equation}
and
	\begin{equation}\label{delta_prod}	
		\Delta_\nu(z)=\prod_{n\geq1}^{}\left(1-\frac{z^2}{\varsigma_{\nu,n}^2}\right),
	\end{equation}
	where $\varsigma_{\nu,n}$ denotes the $n$th positive zero of the function $\Delta_\nu$.
	Taking the logarithmic derivative of both sides of the equation \eqref{delta_prod} we obtain
	\begin{align}\label{delta_ratio1}
		\frac{\Delta_\nu^\prime(z)}{\Delta_\nu(z)}=-2\sum_{k= 0}^{\infty}\rho_{k+1}z^{2k+1},\quad  |z|<\varsigma_{\nu,1},
	\end{align}
	where $\rho_k=\sum_{n\geq1}^{}\varsigma_{\nu,n}^{-2k}$. While, from equation \eqref{delta_series} we have
	\begin{align}\label{delta_ratio2}
		\frac{\Delta_\nu^\prime(z)}{\Delta_\nu(z)}=\left.\sum_{n=0}^{\infty}\xi_nz^{2n+1}\right/\sum_{n=0}^{\infty}\kappa_n z^{2n}
	\end{align}
	for $$\xi_n=\frac{(-1)^{n+1}2(2n+3)^2}{2^{2n+2}n!(\nu+1)_{n+1}}\quad \mbox{and} \quad \kappa_n=\frac{(-1)^n(2n+1)^2}{2^{2n}n!(\nu+1)_n}.$$
	Now, following the result from \cite[Theorem 6]{ABO18} that the radius of convexity $r^c(g_\nu)$ is the smallest positive zero $\varsigma_{\nu,1}$ of the function $z\mapsto\Delta_\nu(z)=(zg_\nu^\prime(z))^\prime$ and by using the Euler-Rayleigh inequalities (cf. \cite{IM95}) we obtain
	\begin{equation}\label{upr_lvr_bnd}
		\rho_k^{-\frac{1}{k}}<\varsigma_{\nu,1}^2=(r^c(g_\nu))^2<\frac{\rho_k}{\rho_{k+1}}\quad \mbox{for all } \nu>-1 \text{ and } k \in \mathbb{N}.
	\end{equation}
	From equations \eqref{delta_ratio1} and \eqref{delta_ratio2} we have that
	\begin{equation}\label{rho_series}
		-2\sum_{k= 0}^{\infty}\rho_{k+1}z^{2k+1}=\left.\sum_{n=0}^{\infty}\xi_nz^{2n+1}\right/\sum_{n=0}^{\infty}\kappa_nz^{2n}.
	\end{equation}
	We express the generalized formula for $\rho_n$ by using potential polynomials. First consider
	\begin{align*}
		\left[\sum_{n\geq 0}^{}\kappa_nz^{2n}\right]^{-1}=\left[1+\sum_{n\geq 1}^{}\kappa_nz^{2n}\right]^{-1}=\sum_{m=0}^{\infty}(-1)^m\left(\sum_{n=1}^{\infty}\kappa_nz^{2n}\right)^m,
	\end{align*}
	where we used the fact that $\kappa_0=1$. Furthermore, we have that
	\begin{align*}
		\left(\sum_{n=1}^{\infty}\kappa_n z^{2n}\right)^m&=\kappa_1^m z^{2m}\left(1+\frac{\kappa_2}{\kappa_1}z^2+\frac{\kappa_3}{\kappa_1}(z^2)^2+\ldots\frac{\kappa_n}{\kappa_1}(z^2)^{n-1}+\ldots\right)^m\\
		&=\kappa_1^m z^{2m}\left(1+\sum_{n=1}^{\infty}f_n z^{2n}\right)^m\\
		&=\kappa_1^m z^{2m}\sum_{n=0}^{\infty}A_{m,n}(f_1,\ldots,f_n)z^{2n},
	\end{align*}
	where
	\begin{equation}\label{f_and_A}
		f_n=\frac{\kappa_{n+1}}{\kappa_1}\quad \text{ and }\quad A_{m,n}(f_1,\ldots,f_n)=\sum_{}^{}\binom{m}{p}\frac{p!}{p_1!\ldots p_n!}f_1^{p_1}f_2^{p_2}\ldots f_n^{p_n},
	\end{equation}
	and the sum extends over all sequences $p_1,\ldots,p_n$ of non-negative integer such that $p_1+2p_2+\ldots+np_n=n$ and $p_1+p_2+\ldots+p_n=p$ (see \cite[Appendix]{Ne13}). Consequently, we write
	\begin{align*}
		\left[\sum_{n\geq 0}^{}\kappa_nz^{2n}\right]^{-1}&=1+\sum_{m=1}^{\infty}(-1)^m\left(\sum_{n=1}^{\infty}\kappa_nz^{2n}\right)^m\\
		&=1+\sum_{m=1}^{\infty}(-1)^m\kappa_1^m z^{2m}\sum_{n=0}^{\infty}A_{m,n}(f_1,\ldots,f_n)z^{2n}\\
		&=1-\kappa_1A_{1,0}z^2+\left(-\kappa_1A_{1,1}(f_1)+\kappa_1^2A_{2,0}\right)z^4+z^{6}\left[\sum_{p=1}^{3}(-1)^p\kappa_1^pA_{3,3-p}(f_1,\ldots,f_{3-p})\right]+\ldots\\
		&=1+\sum_{m=1}^{\infty}\left[\sum_{p=1}^{m}(-1)^p\kappa_1^pA_{m,m-p}(f_1,\ldots,f_{m-p})\right]z^{2m}\\
		&=\sum_{m=0}^{\infty}\pi_m z^{2m},
	\end{align*}
	where
	$\pi_0=1$ and
	\begin{equation}\label{pi_exp}
		\pi_m=\sum_{p=1}^{m}(-1)^p\kappa_1^pA_{m,m-p}(f_1,\ldots,f_{m-p}) \text{ for all }m\in \mathbb{N}.
	\end{equation}
	From equation \eqref{rho_series} we obtain that
	\begin{align*}
		-2\sum_{n=0}^{\infty}\rho_{n+1}z^{2n+1}&=\left.\sum_{n= 0}^{\infty}\xi_nz^{2n+1}\right/\sum_{n=0}^{\infty}\kappa_nz^{2n}\\
		&=\left(\sum_{n= 0}^{\infty}\xi_nz^{2n+1}\right)\left(\sum_{m=0}^{\infty}\pi_m z^{2m}\right)\\
		&=z\left(\sum_{n=0}^{\infty}\sum_{m=0}^{n}\pi_m\xi_{n-m}z^{2n}\right),
	\end{align*}
	which implies that
	\begin{equation}\label{bound_exprsn}
		-2\rho_{n+1}=\sum_{m=0}^{n}\pi_m\xi_{n-m}.
	\end{equation}
	Consequently, we write the upper bound for the radius of convexity of Bessel functions as
	\begin{equation}\label{Upr_bnd}
		\frac{\rho_{n+1}}{\rho_{n+2}}=\left.\sum_{m=0}^{n}\pi_m\xi_{n-m}\right/\sum_{m=0}^{n+1}\pi_m\xi_{n-m+1}.
	\end{equation}
	Since we are going to consider the asymptotic behavior of the radius of convexity of $g_\nu(z)$, as $\nu\rightarrow \infty$, we observe that the next expansions are valid as $\nu\rightarrow \infty$
$$\xi_n=\frac{(-1)^{n+1}2(2n+3)^2}{2^{2n+2}n!(\nu+1)_{n+1}}=\frac{(-1)^{n+1}2(2n+3)^2}{2^{2n+2}n!(\nu+1)(\nu+2)\ldots(\nu+1+n)}=
\frac{1}{\nu^{n+1}}\sum_{m=0}^{\infty}\frac{\xi_m^{(n)}}{\nu^m}$$
and
$$\kappa_n=\frac{(-1)^n(2n+1)^2}{2^{2n}n!(\nu+1)_n}=\frac{(-1)^n(2n+1)^2}{2^{2n}n!(\nu+1)(\nu+2)\ldots(\nu+n)}
		=\frac{1}{\nu^{n}}\sum_{m=0}^{\infty}\frac{\kappa_m^{(n)}}{\nu^m}$$
for some coefficients $\xi_m^{(n)}$ and $\kappa_m^{(n)}$. Moreover, by using the infinite series expansion of $\kappa_n$ and equation \eqref{f_and_A}, we arrive at
$$f_n=\frac{\kappa_{n+1}}{\kappa_1}=\left.\frac{1}{\nu^{n+1}}\sum_{m=0}^{\infty}\frac{\kappa_m^{(n+1)}}{\nu^m}\right/\frac{1}{\nu}\sum_{m=0}^{\infty}\frac{\kappa_m^{(1)}}{\nu^m}=
\frac{1}{\nu^n}\sum_{m=0}^{\infty}\frac{f_m^{(n)}}{\nu^m}$$
for some coefficients $f_m^{(n)}$. In view of the relation $p_1+2p_2+\ldots+np_n=n$ and equation \eqref{f_and_A} we obtain
	\begin{align*}
		A_{m,n}(f_1,\ldots,f_n)&=\sum_{}^{}\binom{m}{p}\frac{p!}{p_1!\ldots p_n!}\left(\frac{1}{\nu}\sum_{m=0}^{\infty}\frac{f_m^{(1)}}{\nu^m}\right)^{p_1}\ldots
		\left(\frac{1}{\nu^n}\sum_{m=0}^{\infty}\frac{f_m^{(n)}}{\nu^m}\right)^{p_n}\\
		&=\frac{1}{\nu^n}\sum_{}^{}\binom{m}{p}\frac{p!}{p_1!\ldots p_n!}\left(\sum_{m=0}^{\infty}\frac{f_m^{(1)}}{\nu^m}\right)^{p_1}\ldots \left(\sum_{m=0}^{\infty}\frac{f_m^{(n)}}{\nu^m}\right)^{p_n}\\
		&=\frac{1}{\nu^n}\sum_{l=0}^{\infty}\frac{A^{m,n}_l}{\nu^l},
	\end{align*}
	for some coefficients $A_l^{(m,n)}$. Now, by using the equation \eqref{pi_exp} and series expansion of $A_{m,n}(f_1,\ldots,f_n)$ we obtain
	\begin{align*}
		\pi_m&=\sum_{p=1}^{m}(-1)^p\kappa_1^pA_{m,m-p}(f_1,\ldots,f_{m-p})\\
		&=\sum_{p=1}^{m}(-1)^p\left(\frac{1}{\nu}\sum_{n=0}^{\infty}\frac{\kappa_n^{(1)}}{\nu^n}\right)^p\left(\frac{1}{\nu^{m-p}}\sum_{l=0}^{\infty}\frac{A^{(m,m-p)}_l}{\nu^l}\right)\\
		&=\frac{1}{\nu^m}\sum_{p=1}^{m}(-1)^p\left(\sum_{n=0}^{\infty}\frac{\kappa_n^{(1)}}{\nu^n}\right)^p\left(\sum_{l=0}^{\infty}\frac{A^{(m,m-p)}_l}{\nu^l}\right)\\
		&=\frac{1}{\nu^m}\sum_{n=0}^{\infty}\frac{\pi_n^{(m)}}{\nu^n},
	\end{align*}
	for some coefficients $\pi_n^{(m)}$. Moreover, in view of equation \eqref{bound_exprsn}, we obtain that
	$$\rho_{n+1}=-\frac{1}{2}\sum_{m=0}^{n}\pi_m\xi_{n-m}=-\frac{1}{2}\sum_{m=0}^{n}\left(\frac{1}{\nu^m}\sum_{n=0}^{\infty}\frac{\pi_n^{(m)}}{\nu^n}\right)
	\left(\frac{1}{\nu^{n-m+1}}\sum_{l=0}^{\infty}\frac{\xi_l^{(n-m)}}{\nu^l}\right)=\frac{1}{\nu^{n+1}}\sum_{l=0}^{\infty}\frac{\rho_l^{(n+1)}}{\nu^l}$$
	for some coefficients $\rho_l^{(n+1)}$.
	Finally, by using \eqref{Upr_bnd} we have that
	\begin{align*}
		\frac{\rho_{n+1}}{\rho_{n+2}}&=\left.\sum_{m=0}^{n}\pi_m\xi_{n-m}\right/\sum_{m=0}^{n+1}\pi_m\xi_{n-m+1}\\
		&=\left.\frac{1}{\nu^{n+1}}\sum_{l=0}^{\infty}\frac{\rho_l^{(n+1)}}{\nu^l}\right/\frac{1}{\nu^{n+2}}\sum_{l=0}^{\infty}\frac{\rho_l^{(n+2)}}{\nu^l}\\
		&=\nu\left(\sum_{l=0}^{\infty}\frac{\rho_l^{(n+1)}}{\nu^l}\right)\left(\sum_{l=0}^{\infty}\frac{\rho_l^{(n+2)}}{\nu^l}\right)^{-1}.
	\end{align*}
	By expanding the series in above equation and multiplying we obtain
	\begin{equation}\label{bnd_form}
		\frac{\rho_{n+1}}{\rho_{n+2}}=\nu\left(d+\mathcal{O}\left(\frac{1}{\nu}\right)\right)
	\end{equation}
	for large $\nu$ and some constant $d$. Moreover, the left-hand side of inequality \eqref{upr_lvr_bnd} also admits a similar asymptotic form as $\nu\to \infty$ as given below
$$\rho_{n+1}^{-\frac{1}{n+1}}=
\left(\frac{1}{\nu^{n+1}}\sum_{l=0}^{\infty}\frac{\rho_l^{(n+1)}}{\nu^l}\right)^{-\frac{1}{n+1}}
=\nu\cdot\left(\sum_{l=0}^{\infty}\frac{\rho_l^{(n+1)}}{\nu^l}\right)^{-\frac{1}{n+1}}$$
	or we can write
	\begin{equation}\label{asym_lwr_bnd}
		\rho_{n+1}^{-\frac{1}{n+1}}=\nu\left(e+\mathcal{O}\left(\frac{1}{\nu}\right)\right),
	\end{equation}
	for large $\nu$ and some constant $e$.
	Note that the asymptotic form of the ratio $\frac{\rho_{n+1}}{\rho_{n+2}}$ and $\rho_{n+1}^{-\frac{1}{n+1}}$ are true for all $n\in \mathbb{N}$. Moreover, from \cite[Lemma 3.2]{IM95} the left-hand side of the inequality \eqref{upr_lvr_bnd} increases and the right-hand side ratio of \eqref{upr_lvr_bnd} decreases to $(r^c(g_\nu))^2$ as $n\to \infty$. Also these bounds will admit the asymptotic form \eqref{bnd_form} and \eqref{asym_lwr_bnd} for large $\nu$. From equations \eqref{bnd_form}, \eqref{asym_lwr_bnd} for large $\nu$ and $n\to \infty$,  we conclude that radius of convexity  $r^c(g_\nu)$ admits the asymptotic form as
	\begin{align*}
		\left(r^c(g_\nu)\right)^2=\nu\left(c+\mathcal{O}\left(\frac{1}{\nu}\right)\right).
	\end{align*}
	for some positive constant $c$, since the radius of convexity cannot be negative.
\end{proof}

\begin{proof}[\bf Proof of Lemma \ref{Lemma4}]
We can prove this Lemma by using \cite[Theorem 7]{ABO18}. The proof is very similar to the proof of the previous Lemma, so we omit the details.
\end{proof}

\begin{proof}[\bf Proof of Lemma \ref{unif_bounds_g}]
	From the proof of \cite[Theorem 3.2. i.]{DS17} for $\nu>-1$, the radius of uniform convexity of the function $g_\nu(z)$ is the smallest positive zero of the equation
	\begin{align*}
		1+2r\frac{(2\nu -1)J_{\nu+1}(r)-rJ_\nu(r)}{J_\nu(r)-rJ_{\nu+1}(r)}=0.
	\end{align*}
	Also, by using the recurrence relation
	%	$2\nu J_\nu(z)=z[J_{\nu-1}(z)+J_{\nu+1}(z)]$
	\begin{equation}\label{unif_bessel_reln}
		2\nu J_\nu(z)=z[J_{\nu-1}(z)+J_{\nu+1}(z)]
	\end{equation}
	and definition of $g_\nu(z)$, we write
	\begin{equation}\label{unif_g_nu_expres}
		1+2\frac{rg^{\prime\prime}_\nu(r)}{g_\nu^\prime(r)}=1+2r\frac{(2\nu -1)J_{\nu+1}(r)-rJ_\nu(r)}{J_\nu(r)-rJ_{\nu+1}(r)}.
	\end{equation}
	%	 function $1+2\frac{rg^{\prime\prime}_\nu(r)}{g_\nu^\prime(r)}$.
	From the above discussion we conclude that for $\nu>-1$ the radius of uniform convexity of $g_\nu(z)$ is the smallest positive root of equation
	$$
	g_\nu^\prime(z)+2zg_\nu^{\prime\prime}(z)=0.
	$$
	By using the series form of Bessel functions \eqref{Bessel_series}, its derivative \eqref{der_Bessel_series} and the definition of $g_\nu(z)$, we obtain
	\begin{equation}\label{unif_Delta_series}
		\Phi_\nu(z):=g_\nu^\prime(z)+2zg_\nu^{\prime\prime}(z)=1+\sum_{n=1}^{\infty}\frac{(-1)^n(2n+1)(4n+1)}{2^{2n}n!(\nu+1)_n}z^{2n}.
	\end{equation}

Since $g_\nu(z)$ is a member of $\mathcal{LP}$, which is closed under differentiation, it follows that  $2zg_\nu^\prime(z)$ and therefore $2g_\nu^\prime(z)+2zg_\nu^{\prime\prime}(z)$ also belongs to the class $\mathcal{LP}.$ Moreover, $2zg_\nu'(z)$ satisfies all the conditions of the Laguerre separation theorem \cite[Theorem 2.8.1]{Bo54}, consequently the zeros of $2g_\nu^\prime(z)+2zg_\nu^{\prime\prime}(z)$ are separated by the zeros of $2zg_\nu^\prime(z)$ as well as $g_\nu^\prime(z)$. Therefore, we conclude that zeros of $g_\nu^\prime(z)+2zg_\nu^{\prime\prime}(z)$ are all real. Furthermore, the larger zeros of $g_\nu^\prime(z)+2zg_\nu^{\prime\prime}(z)$ correspond to larger argument $z$ and from \cite[eq.1]{Sk02} we can conclude that $g_\nu^\prime(z)+2zg_\nu^{\prime\prime}(z)$ has growth similar to $\cos z$. By the above discussion we conclude that the function $\Phi_\nu(z)$ also belongs to the Laguerre-P\'olya class. Another approach is to show via recurrence relations and the Bessel differential equation that the equation $g_\nu^\prime(z)+2zg_\nu^{\prime\prime}(z)=0$ is equivalent to $(4\nu-3)zJ_{\nu+1}(z)=(2z^2-1)J_{\nu}(z)$ and in view of the well-known Mittag-Leffler expansion for $J_{\nu+1}(z)/J_{\nu}(z)$ the above equation has only real solutions.

Now, let $\gamma_{\nu,n}$ represents the $n$th positive zero of $\Phi_\nu(z)$. Then, the function $\Phi_\nu(z)$ can be expressed as follows
	\begin{equation}\label{unif_Delta_inf_prod}
		\Phi_\nu(z)=\prod_{n\geq 1}^{}\left(1-\frac{z^2}{\gamma_{\nu,n}^2}\right).
	\end{equation}
	By taking the logarithmic derivative of \eqref{unif_Delta_inf_prod} we obtain that
	\begin{equation}\label{unif_log_der_delta}
		\frac{\Phi_\nu^\prime(z)}{\Phi_\nu(z)}=-2\sum_{k\geq 0}^{}\omega_{k+1}z^{2k+1},\quad  |z|<\gamma_{\nu,n}^2,
	\end{equation}
	where $\omega_{k}=\sum_{n\geq 1}^{}\gamma_{\nu,n}^{-2k}$. Also, by considering the infinite sum representation \eqref{unif_Delta_series} we obtain that
	\begin{equation}\label{unif_div_delta}
		\frac{\Phi_\nu^\prime(z)}{\Phi_\nu(z)}=\left.\sum_{n\geq 0}^{}\xi_nz^{2n+1}\right/\sum_{n\geq 0}^{}\kappa_nz^{2n},
	\end{equation}
	where
	\begin{equation}\label{unif_xi}
		\xi_n= \frac{(-1)^{n+1}2(2n+3)(4n+5)}{2^{2n+2}n!(\nu+1)_{n+1}}
	\end{equation}
	and
	\begin{equation}\label{unif_kappa}
		\kappa_n= \frac{(-1)^{n}(2n+1)(4n+1)}{2^{2n}n!(\nu+1)_{n}}.
	\end{equation}
	By using the Euler-Rayleigh inequalities (cf. \cite{IM95}), equations \eqref{unif_log_der_delta} and \eqref{unif_div_delta} for $\nu>-1$ and $k \in \mathbb{N}$ we obtain that
	\begin{equation}\label{unif_upr_lvr_bnd}
		\omega_k^{-\frac{1}{k}}<\gamma_{\nu,1}^2=(r^{uc}(g_\nu))^2<\frac{\omega_k}{\omega_{k+1}}
	\end{equation}
	and
	\begin{equation}\label{unif_rho_series}
		-2\sum_{k= 0}^{\infty}\omega_{k+1}z^{2k+1}=\left.\sum_{n=0}^{\infty}\xi_nz^{2n+1}\right/\sum_{n=0}^{\infty}\kappa_nz^{2n}.
	\end{equation}
	
Now, we find a generalized formula for $\omega_n$ by using potential polynomials. First consider
$$\left[\sum_{n\geq 0}^{}\kappa_nz^{2n}\right]^{-1}=\left[1+\sum_{n\geq 1}^{}\kappa_nz^{2n}\right]^{-1}=\sum_{m=0}^{\infty}(-1)^m\left(\sum_{n=1}^{\infty}\kappa_nz^{2n}\right)^m,$$
where we used the fact that $\kappa_0=1$. In the rest of the proof, we use a similar expansion technique as in the proof of Lemma \ref{Lemma_g_asy}. We observe that
\begin{align*}
\left(\sum_{n=1}^{\infty}\kappa_n z^{2n}\right)^m=\kappa_1^m z^{2m}\sum_{n=0}^{\infty}A_{m,n}(f_1,\ldots,f_n)z^{2n},
\end{align*}
	for
	\begin{equation}\label{unif_f_and_A}
		f_n=\frac{\kappa_{n+1}}{\kappa_1}\quad \text{ and }\quad A_{m,n}(f_1,\ldots,f_n)=\sum_{}^{}\binom{m}{p}\frac{p!}{p_1!\ldots p_n!}f_1^{p_1}f_2^{p_2}\ldots f_n^{p_n},
	\end{equation}
	where the sum extends over all sequences $p_1,\ldots,p_n$ of non-negative integer such that $p_1+2p_2+\ldots+np_n=n$ and $p_1+p_2+\ldots+p_n=p$ (see \cite[Appendix]{Ne13}). Consequently, we arrive at
$$\left[\sum_{n\geq 0}^{}\kappa_nz^{2n}\right]^{-1}=1+\sum_{m=1}^{\infty}(-1)^m\left(\sum_{n=1}^{\infty}\kappa_nz^{2n}\right)^m=\sum_{m=0}^{\infty}\pi_m z^{2m},$$
	where
	$\pi_0=1$ and
	\begin{equation}\label{unif_pi_exp}
		\pi_m=\sum_{p=1}^{m}(-1)^p\kappa_1^pA_{m,m-p}(f_1,\ldots,f_{m-p}) \text{ for all }m\in \mathbb{N}.
	\end{equation}
	From equation \eqref{unif_rho_series} we obtain that 	
	\begin{align*}
		-2\sum_{n=0}^{\infty}\omega_{n+1}z^{2n+1}&=\left.\sum_{n= 0}^{\infty}\xi_nz^{2n+1}\right/\sum_{n=0}^{\infty}\kappa_nz^{2n}\\
		&=\left(\sum_{n= 0}^{\infty}\xi_nz^{2n+1}\right)\left(\sum_{m=0}^{\infty}\pi_m z^{2m}\right)\\
		&=z\left(\sum_{n=0}^{\infty}\sum_{m=0}^{n}\pi_m\xi_{n-m}z^{2n}\right)
	\end{align*}
	\begin{equation}\label{unif_bound_exprsn}
		-2\omega_{n+1}=\sum_{m=0}^{n}\pi_m\xi_{n-m}.
	\end{equation}
Now, for $k\in \mathbb{N}$, by considering the Euler-Rayleigh inequalities $\omega_{k}^{-\frac{1}{k}}<\gamma_{\nu,1}^2<\frac{\omega_{k}}{\omega_{k+1}}$ for $\nu>-1$ we obtain the required bounds for the radius of uniform convexity of $g_\nu(z)$. In particular, for $k=1$ we obtain the inequality \eqref{unif_first_inq}.
\end{proof}

\begin{proof}[\bf Proof of Lemma \ref{unif_bounds_h}]
	From the proof of \cite[Theorem 3.3. i.]{DS17} for $\nu>-1$, the radius of uniform convexity of the function $h_\nu(z)$ is the smallest positive zero of the equation
	\begin{align*}
		1+r^{\frac{1}{2}}\frac{2(\nu -1)J_{\nu+1}(r^{\frac{1}{2}})-r^{\frac{1}{2}}J_\nu(r^{\frac{1}{2}})}{2J_\nu(r^{\frac{1}{2}})-r^{\frac{1}{2}}J_{\nu+1}(r^{\frac{1}{2}})}=0.
	\end{align*}
	Also, by using the recurrence relation \eqref{unif_bessel_reln} and the definition of $h_\nu(z)$, we write
	\begin{equation}\label{unif_h_nu_expres}
		1+2\frac{rh^{\prime\prime}_\nu(r)}{h_\nu^\prime(r)}=1+r^{\frac{1}{2}}\frac{2(\nu -1)J_{\nu+1}(r^{\frac{1}{2}})-r^{\frac{1}{2}}J_\nu(r^{\frac{1}{2}})}{2J_\nu(r^{\frac{1}{2}})-r^{\frac{1}{2}}J_{\nu+1}(r^{\frac{1}{2}})}.
	\end{equation}
	From the above discussion we conclude that for $\nu>-1$ the radius of uniform convexity of $h_\nu(z)$ is the smallest positive root of the equation
	$$
	h_\nu^\prime(z)+2zh_\nu^{\prime\prime}(z)=0.
	$$
	Now, in view of the infinite series representations of the Bessel function \eqref{Bessel_series}, its derivative \eqref{der_Bessel_series} and the definition of $h_\nu(z)$, we obtain
	\begin{equation}\label{unif_Delta_series_h}
		\Theta_\nu(z):=h_\nu^\prime(z)+2zh_\nu^{\prime\prime}(z)=1+\sum_{n=1}^{\infty}\frac{(-1)^n(n+1)(2n+1)}{2^{2n}n!(\nu+1)_n}z^n.
	\end{equation}
Since the function $h_\nu$ is a member of the Laguerre-P\'olya class (denoted as $\mathcal{LP}$) of entire functions and the class is closed under differentiation, by using a similar argument as for \eqref{unif_Delta_series}, we conclude that the function $\Theta_\nu(z)$ will be also in the $\mathcal{LP}$ class. Consequently, all zeros of the function $\Theta_\nu(z)$ are real. Let $\delta_{\nu,n}$ represents the $n$th positive zero of $\Theta_\nu(z)$. The function $\Theta_\nu(z)$ can be expressed through an infinite product as follows
	\begin{equation}\label{unif_Delta_inf_prod_h}
		\Theta_\nu(z)=\prod_{n\geq 1}^{}\left(1-\frac{z}{\delta_{\nu,n}}\right).
	\end{equation}
	By taking the logarithmic derivative of \eqref{unif_Delta_inf_prod_h} we arrive at
	\begin{equation}\label{unif_log_der_delta_h}
		\frac{\Theta_\nu^\prime(z)}{\Theta_\nu(z)}=-\sum_{k\geq 0}^{}\sigma_{k+1}z^k,\quad |z|<\delta_{\nu,1},
	\end{equation}
	where $\sigma_{k}=\sum_{n\geq 1}^{}\delta_{\nu,n}^{-k}$. Also, by considering the infinite sum representation in \eqref{unif_Delta_series_h} we obtain that
	\begin{equation}\label{unif_div_delta_h}
		\frac{\Theta_\nu^\prime(z)}{\Theta_\nu(z)}=\left.\sum_{n\geq 0}^{}\lambda_n z^{n}\right/\sum_{n\geq 0}^{}\mu_nz^{n},
	\end{equation}
	where
	\begin{align*}
		\lambda_n= \frac{(-1)^{n+1}(n+1)(n+2)(2n+3)}{2^{2n+2}(n+1)!(\nu+1)_{n+1}}\quad \text{ and }\quad \mu_n= \frac{(-1)^n(n+1)(2n+1)}{2^{2n}n!(\nu+1)_n}.
	\end{align*}
	By using the Euler-Rayleigh inequalities (cf. \cite{IM95}), equations \eqref{unif_log_der_delta_h} and \eqref{unif_div_delta_h}, for all $\nu>-1$ and $k \in \mathbb{N}$ we obtain that
	\begin{align*}
		\sigma_k^{-\frac{1}{k}}<\delta_{\nu,1}=r^{uc}(h_\nu)<\frac{\sigma_k}{\sigma_{k+1}}
	\end{align*}
	and
	\begin{equation}\label{unif_rho_series_h}
		-\sum_{k= 0}^{\infty}\sigma_{k+1}z^{k}=\left.\sum_{n=0}^{\infty}\lambda_nz^{n}\right/\sum_{n=0}^{\infty}\mu_nz^{n}.
	\end{equation}
	Next, we find a generalized formula for $\sigma_n$ by using potential polynomials. First consider
	\begin{align*}
		\left[\sum_{n\geq 0}^{}\mu_nz^{n}\right]^{-1}=\left[1+\sum_{n\geq 1}^{}\mu_nz^{n}\right]^{-1}=\sum_{m=0}^{\infty}(-1)^m\left(\sum_{n=1}^{\infty}\mu_nz^{n}\right)^m,
	\end{align*}
	where we used the fact that $\mu_0=1$. Moreover, observe that
	\begin{align*}
		\left(\sum_{n=1}^{\infty}\mu_n z^{n}\right)^m&=\mu_1^m z^{m}\left(1+\frac{\mu_2}{\mu_1}z+\frac{\mu_3}{\mu_1}z^2+\ldots\frac{\mu_n}{\mu_1}z^{n-1}+\ldots\right)^m\\
		&=\mu_1^m z^{m}\left(1+\sum_{n=1}^{\infty}f_nz^n\right)^m\\
		&=\mu_1^m z^{m}\sum_{n=0}^{\infty}A_{m,n}(f_1,\ldots,f_n)z^n,
	\end{align*}
	for
	\begin{align*}
		f_n=\frac{\mu_{n+1}}{\mu_1}\quad \text{ and }\quad A_{m,n}(f_1,\ldots,f_n)=\sum_{}^{}\binom{m}{p}\frac{p!}{p_1!\ldots p_n!}f_1^{p_1}f_2^{p_2}\ldots f_n^{p_n},
	\end{align*}
	where the sum extends over all sequences $p_1,\ldots,p_n$ of non-negative integer such that $p_1+2p_2+\ldots+np_n=n$ and $p_1+p_2+\ldots+p_n=p$ (see \cite[Appendix]{Ne13}). Similarly, as in the proof of Lemma \ref{unif_bounds_g}, we can write
	\begin{align*}
		\left[\sum_{n\geq 0}^{}\mu_nz^{n}\right]^{-1}=1+\sum_{m=1}^{\infty}(-1)^m\left(\sum_{n=1}^{\infty}\mu_nz^{n}\right)^m=\sum_{m=0}^{\infty}\pi_m z^{m},
	\end{align*}
	where
	$\pi_0=1$ and
	\begin{align*}
		\pi_m=\sum_{p=1}^{m}(-1)^p\mu_1^pA_{m,m-p}(f_1,\ldots,f_{m-p})\quad  \text{ for all }m\in \mathbb{N}.
	\end{align*}
	From equation \eqref{unif_rho_series_h} we obtain
	\begin{align*}
		-\sum_{n=0}^{\infty}\sigma_{n+1}z^{n}&=\left.\sum_{n= 0}^{\infty}\lambda_nz^n\right/\sum_{n=0}^{\infty}\mu_nz^n=\sum_{n=0}^{\infty}\sum_{m=0}^{n}\pi_m\lambda_{n-m}z^{n}
	\end{align*}
	and thus
	\begin{equation}\label{unif_bound_exprsn_h}
		-\sigma_{n+1}=\sum_{m=0}^{n}\pi_m\lambda_{n-m}.
	\end{equation}
Now, for $k\in \mathbb{N}$, by considering the Euler-Rayleigh inequalities $\sigma_{k}^{-\frac{1}{k}}<\delta_{\nu,1}<\frac{\sigma_{k}}{\sigma_{k+1}}$ for $\nu>-1$ we arrive to the required bounds for the radius of uniform convexity of $h_\nu(z)$.	
\end{proof}

\begin{proof}[\bf Proof of Lemma \ref{unif_Lemma_g_asy}]
	From Lemma \ref{unif_bounds_g} and equation \eqref{unif_bound_exprsn} we can write the upper bound for the radius of uniform convexity of $g_\nu(z)$ as
	\begin{equation}\label{unif_Upr_bnd}
		\frac{\omega_{n+1}}{\omega_{n+2}}=\left.\sum_{m=0}^{n}\pi_m\xi_{n-m}\right/\sum_{m=0}^{n+1}\pi_m\xi_{n-m+1}.
	\end{equation}
	Since we want to discuss the asymptotic behavior of the radius of uniform convexity of $g_\nu(z)$, we are considering $\nu$ large enough so that the expansions below are valid. Notice that from equations \eqref{unif_xi} and \eqref{unif_kappa} we obtain
	\begin{equation}\label{unif_xi_series}
		\xi_n=\frac{(-1)^{n}(2n+1)(4n+1)}{2^{2n}n!(\nu+1)_{n}}=\frac{(-1)^{n}(2n+1)(4n+1)}{2^{2n}n!(\nu+1)(\nu+2)\ldots(\nu+1+n)}=\frac{1}{\nu^{n+1}}\sum_{m=0}^{\infty}\frac{\xi_m^{(n)}}{\nu^m}
	\end{equation}
	and
	\begin{equation}\label{unif_kappa_series}
		\kappa_n=\frac{(-1)^{n+1}2(2n+3)(4n+5)}{2^{2n+2}n!(\nu+1)_{n+1}}=\frac{(-1)^{n+1}2(2n+3)(4n+5)}{2^{2n+2}n!(\nu+1)(\nu+2)\ldots(\nu+n)}=\frac{1}{\nu^{n}}\sum_{m=0}^{\infty}\frac{\kappa_m^{(n)}}{\nu^m},
	\end{equation}
	for some coefficients $\xi_m^{(n)}$ and $\kappa_m^{(n)}$. By using the infinite series expansion of $\kappa_n$ \eqref{unif_kappa_series} and equation  \eqref{unif_f_and_A}
	\begin{align*}
		f_n=\frac{\kappa_{n+1}}{\kappa_1}=\left.\frac{1}{\nu^{n+1}}\sum_{m=0}^{\infty}\frac{\kappa_m^{(n+1)}}{\nu^m}\right/\frac{1}{\nu}\sum_{m=0}^{\infty}\frac{\kappa_m^{(1)}}{\nu^m}=\frac{1}{\nu^n}\sum_{m=0}^{\infty}\frac{f_m^{(n)}}{\nu^m},
	\end{align*}
	for some coefficients $f_m^{(n)}$. In the rest of the proof, we use a similar expansion technique as in the proof of Lemma \ref{Lemma_g_asy}. In view of the relation $p_1+2p_2+\ldots+np_n=n$  and equation  \eqref{unif_f_and_A} we obtain that
\begin{align*}
A_{m,n}(f_1,\ldots,f_n)&=\sum_{}^{}\binom{m}{p}\frac{p!}{p_1!\ldots p_n!}\left(\frac{1}{\nu}\sum_{m=0}^{\infty}\frac{f_m^{(1)}}{\nu^m}\right)^{p_1}\ldots \left(\frac{1}{\nu^n}\sum_{m=0}^{\infty}\frac{f_m^{(n)}}{\nu^m}\right)^{p_n}=\frac{1}{\nu^n}\sum_{l=0}^{\infty}\frac{A^{m,n}_l}{\nu^l},
\end{align*}
for some coefficients $A_l^{(m,n)}$. Now by using the equation \eqref{unif_pi_exp} and series expansion of $A_{m,n}(f_1,\ldots,f_n)$ we obtain that
\begin{align*}
\pi_m&=\sum_{p=1}^{m}(-1)^p\kappa_1^pA_{m,m-p}(f_1,\ldots,f_{m-p})=\frac{1}{\nu^m}\sum_{n=0}^{\infty}\frac{\pi_n^{(m)}}{\nu^n},
\end{align*}
for some coefficients $\pi_n^{(m)}$. Moreover, from equation \eqref{unif_bound_exprsn}, the above expansion of $\pi_m$ and the expression \eqref{unif_xi_series} of $\xi_n$, we arrive at
	$$\omega_{n+1}=-\frac{1}{2}\sum_{m=0}^{n}\pi_m\xi_{n-m}=-\frac{1}{2}\sum_{m=0}^{n}\left(\frac{1}{\nu^m}
\sum_{n=0}^{\infty}\frac{\pi_n^{(m)}}{\nu^n}\right)\left(\frac{1}{\nu^{n-m+1}}\sum_{l=0}^{\infty}\frac{\xi_l^{(n-m)}}{\nu^l}\right)
=\frac{1}{\nu^{n+1}}\sum_{l=0}^{\infty}\frac{\omega_l^{(n+1)}}{\nu^l},$$
	for some coefficients $\omega_l^{(n+1)}$.
	Finally, in view of \eqref{unif_Upr_bnd} we obtain that
		\begin{align*}
		\frac{\omega_{n+1}}{\omega_{n+2}}=\left.\sum_{m=0}^{n}\pi_m\xi_{n-m}\right/\sum_{m=0}^{n+1}\pi_m\xi_{n-m+1}=
\nu\left(\sum_{l=0}^{\infty}\frac{\omega_l^{(n+1)}}{\nu^l}\right)\left(\sum_{l=0}^{\infty}\frac{\omega_l^{(n+2)}}{\nu^l}\right)^{-1}.
	\end{align*}
	By expanding the series in above equation and multiplying we obtain that
	\begin{equation}\label{unif_bnd_form}
		\frac{\omega_{n+1}}{\omega_{n+2}}=\nu\left(d+\mathcal{O}\left(\frac{1}{\nu}\right)\right)
	\end{equation}
	for large $\nu$ and some constant $d$. Moreover, the left-hand side of inequality \eqref{unif_upr_lvr_bnd} also admit the similar asymptotic form as $\nu\to \infty$ as given below
$$\omega_{n+1}^{-\frac{1}{n+1}}=\left(\frac{1}{\nu^{n+1}}\sum_{l=0}^{\infty}\frac{\omega_l^{(n+1)}}{\nu^l}\right)^{-\frac{1}{n+1}}
=\nu\cdot\left(\sum_{l=0}^{\infty}\frac{\omega_l^{(n+1)}}{\nu^l}\right)^{-\frac{1}{n+1}}$$
or we can write
	\begin{equation}\label{unif_asym_lwr_bnd}
		\omega_{n+1}^{-\frac{1}{n+1}}=\nu\left(e+\mathcal{O}\left(\frac{1}{\nu}\right)\right),
	\end{equation}
for large $\nu$ and some constant $e$. Note that the asymptotic forms of the ratio $\frac{\omega_{n+1}}{\omega_{n+2}}$ and $\omega_{n+1}^{-\frac{1}{n+1}}$ are true for all $n\in \mathbb{N}$. Moreover, from \cite[Lemma 3.2]{IM95} the left-hand side of the inequality \eqref{unif_upr_lvr_bnd} increases and the right-hand side ratio of \eqref{unif_upr_lvr_bnd} decreases to $(r^{uc}(g_\nu))^2$ as $n\to \infty$. Also these bounds will admit the asymptotic form \eqref{unif_bnd_form} and \eqref{unif_asym_lwr_bnd} for large $\nu$. From equations \eqref{unif_bnd_form}, \eqref{unif_asym_lwr_bnd} for large $\nu$ and $n\to \infty$,  we conclude that the radius of uniform convexity  $r^{uc}(g_\nu)$ admits the asymptotic form as
	\begin{align*}
	\left(r^{uc}(g_\nu)\right)^2=\nu\left(c+\mathcal{O}\left(\frac{1}{\nu}\right)\right).
	\end{align*}
	for some positive constant $c$, since the radius of uniform convexity cannot be negative.
\end{proof}

\begin{proof}[\bf Proof of Lemma \ref{unif_Lemma_h_asy}]
The proof of this Lemma is very similar to the proof of Lemma \ref{unif_Lemma_g_asy}, so we omit the details.
\end{proof}

\begin{proof}[\bf Proof of Theorem \ref{Theorem1}]
From the proof of \cite[Lemma 2.4]{BS14} we have
\begin{equation}\label{g_exp}
1+\frac{zg_\nu^{\prime\prime}(z)}{g_\nu^\prime(z)}=1+z\frac{zJ_{\nu+2}(z)-3J_{\nu+1}(z)}{J_\nu(z)-zJ_{\nu+1}(z)}=1-\sum_{n=1}^{\infty}\frac{2z^2}{\alpha_{\nu,n}^2-z^2},
\end{equation}
where $\alpha_{\nu,n}$ is the $n$th positive zero of the Dini function $d_\nu(z)$, defined in \eqref{dini_d}. For $\alpha=0$ \cite[Theorem 1.2]{BS14} implies that $z=r^c(g_\nu)$ is the smallest positive root of the expression in \eqref{g_exp}. Now, since the expression in \eqref{g_exp} is equal to zero at $z=r^c(g_\nu)$ we obtain that
\begin{equation}\label{eqn24}
1-\sum_{n=1}^{\infty}\frac{2(r^c(g_\nu))^2}{\alpha_{\nu,n}^2-(r^c(g_\nu))^2}=0.
\end{equation}
Now with the help of Lemma \ref{Lemma_g_asy} we have
$$\left(r^c(g_\nu)\right)^2=\nu\left(c+\mathcal{O}\left(\frac{1}{\nu}\right)\right)=\nu\left(c+\epsilon(\nu)\right)$$
for large $\nu$, where $c$ is some constant and $\epsilon(\nu)=\mathcal{O}\left(\frac{1}{\nu}\right)$. Rearranging \eqref{eqn24} we find that
\begin{align*}
1&=2\sum_{n\geq1}^{}\frac{(r^c(g_\nu))^2}{\alpha_{\nu,n}^2-(r^c(g_\nu))^2}=2\sum_{n\geq1}^{}\frac{(\nu\left(c+\epsilon(\nu)\right))}{\alpha_{\nu,n}^2-(\nu\left(c+\epsilon(\nu)\right))}\\
&=2\sum_{n\geq1}^{}\frac{1}{\alpha_{\nu,n}^2}\frac{\nu\left(c+\epsilon(\nu)\right)}{1-\frac{\left(c+\epsilon(\nu)\right)}{\alpha_{\nu,n}^2}}\\
&=2\sum_{n\geq1}^{}\frac{\nu\left(c+\epsilon(\nu)\right)}{\alpha_{\nu,n}^2}\sum_{m\geq0}^{}\frac{\left(\nu\left(c+\epsilon(\nu)\right)\right)^m}{\alpha_{\nu,n}^{2m}}\\
&=2\sum_{m\geq0}^{}\left(\nu\left(c+\epsilon(\nu)\right)\right)^{m+1}\sum_{n\geq1}^{}\frac{1}{\alpha_{\nu,n}^{2m+2}}\\
&=2\sum_{m\geq0}^{}\left(\nu\left(c+\epsilon(\nu)\right)\right)^{m+1}\eta_{m+1}(\nu),
\end{align*}
which can be rewritten as
\begin{equation}\label{main_eqn}
\frac{1}{2}=\nu\left(c+\epsilon(\nu)\right)\eta_{1}(\nu)+\sum_{m=1}^{\infty}\left(\nu\left(c+\epsilon(\nu)\right)\right)^{m+1}\eta_{m+1}(\nu),
\end{equation}
provided $\nu$ is sufficiently large. Now, we write
\begin{equation}\label{epsilon_exp}
\epsilon(\nu)=\sum_{n=1}^{N-1}\frac{\epsilon_n}{\nu^n}+R_N(\nu),
\end{equation}
where the coefficients $\epsilon_n$ are given by the recurrence relation \eqref{epsilon_rln}.

An important observation is that given $\epsilon(\nu)=\mathcal{O}\left(\frac{1}{\nu}\right)$, it follows from \eqref{epsilon_exp} that $R_N(\nu)\to 0$ as $\nu\to \infty$. Keeping this in mind, without loss of generality, for fixed positive integer $N$, we assume $R_N(\nu)=\mathcal{O}_N\left(f_N(\nu)\right)$ for some function $f_N(\nu)$. Notice that we can write
\begin{equation}\label{prop_bigO}
	 \mathcal{O}_N\left(f_N(\nu)\right)+\frac{c}{\nu}\mathcal{O}_N\left(f_N(\nu)\right)=\mathcal{O}_N\left(f_N(\nu)\right).
\end{equation}
We shall prove by induction on $N$ that $R_N(\nu)=\mathcal{O}_N\left(\frac{1}{\nu^N}\right)$ for any $N\geq 1$ as $\nu\to \infty$.
Throughout this paper, we use subscripts in the $\mathcal{O}$ notations to indicate the dependence of the implied constant on certain parameters. The statement is true for $N=1$ since $R_1(\nu)=\epsilon(\nu)$=$\mathcal{O}\left(\frac{1}{\nu}\right)$. Let $N\geq2$ and suppose that the statement holds for all $R_k(\nu)$ with $1\leq k\leq N-1$. In view of equation \eqref{epsilon_exp}, Lemma \ref{lemma_eta}, the assumption $R_N(\nu)=\mathcal{O}_N\left(f_N(\nu)\right)$ and the relation (see \cite[p. 2]{BPS14})
\begin{align*}
	\eta_{1}(\nu)=\sum_{m=1}^{\infty}\frac{1}{\alpha_{\nu,n}^2}=\frac{3}{4(\nu+1)},
\end{align*}
the first term on the right-hand side of the equation \eqref{main_eqn} can be expressed as

\begin{align*}
\nu\left(c+\epsilon(\nu)\right)\eta_{1}(\nu)&=\nu\left(c+\sum_{n=1}^{N-1}\frac{\epsilon_n}{\nu^n}+R_N(\nu)\right)\eta_{1}(\nu)\\
&=\nu\left(c+\sum_{n=1}^{N-1}\frac{\epsilon_n}{\nu^n}\right)\eta_{1}(\nu)+\nu R_N(\nu)\eta_{1}(\nu)\\
&=\frac{3}{4}\left(c+\sum_{n=1}^{N-1}\frac{\epsilon_n}{\nu^n}\right)\sum_{n=0}^{\infty}\frac{(-1)^n}{\nu^n}+\nu R_N(\nu)\eta_{1}(\nu)\\
&=\frac{3c}{4}+\frac{3c}{4}\sum_{n=1}^{\infty}\frac{(-1)^n}{\nu^n}+\frac{3}{4\nu}\left(\sum_{n=0}^{N-2}\frac{\epsilon_{n+1}}{\nu^n}\right)\left(\sum_{n=0}^{\infty}\frac{(-1)^n}{\nu^n}\right)+\nu R_N(\nu)\left(\frac{1}{\nu}\sum_{n=0}^{\infty}\frac{\eta_n^{(1)}}{\nu^n}\right)\\
&=\frac{3c}{4}+\frac{3c}{4\nu}\sum_{n=0}^{N-2}\frac{(-1)^{n+1}}{\nu^n}+\frac{3}{4\nu}\sum_{n=0}^{N-2}\sum_{m=0}^{n}\frac{(-1)^{n-m}\epsilon_{m+1}}{\nu^n}+\mathcal{O}_N\left(f_N(\nu)\right)+\mathcal{O}_N\left(\frac{1}{\nu^N}\right)\\
&=\frac{3c}{4}+\frac{3}{4\nu}\sum_{n=0}^{N-2}\left((-1)^{n+1}c+\sum_{m=0}^{n}(-1)^{n-m}\epsilon_{m+1}\right)\frac{1}{\nu^n}+\mathcal{O}_N\left(f_N(\nu)\right)+\mathcal{O}_N\left(\frac{1}{\nu^N}\right).
\end{align*}
With the help of equation \eqref{prop_bigO}, Lemma \ref{lemma_eta} and $R_N(\nu)=\mathcal{O}_N\left(f_N(\nu)\right)$, we simplify the second term in the right-hand side of \eqref{main_eqn} as

\begin{align*}
&\sum_{m=1}^{\infty}\left(\nu(c+\epsilon(\nu))\right)^{m+1}\eta_{m+1}(\nu)=\sum_{m=1}^{\infty}\nu^{m+1}\left(c+\sum_{n=1}^{N-1}\frac{\epsilon_n}{\nu^n}+R_{N}(\nu)\right)^{m+1}\left(\frac{1}{\nu^{m+1}}\sum_{n=0}^{\infty}\frac{\eta_n^{(m+1)}}{\nu^n}\right)\\
&=\sum_{m=1}^{\infty}\left(c+\sum_{n=1}^{N-1}\frac{\epsilon_n}{\nu^n}+R_N(\nu)\right)^{m+1}\left(\sum_{n=0}^{\infty}\frac{\eta_n^{(m+1)}}{\nu^n}\right)\\
&=\sum_{m=1}^{\infty}\left[\left(c+\sum_{n=1}^{N-1}\frac{\epsilon_n}{\nu^n}\right)^{m+1}+ \mathcal{O}_N(f_N(\nu))\right]\left(\sum_{n=0}^{\infty}\frac{\eta_n^{(m+1)}}{\nu^n}\right)\\
&=\sum_{m=1}^{\infty}\left[\left(\sum_{n=0}^{N-1}\frac{A_{m+1,n}\left(\epsilon_1,\ldots,\epsilon_n\right)}{\nu^n}\right)\left(\sum_{k=0}^{\infty}\frac{\eta_k^{(m+1)}}{\nu^k}\right)\right]+\mathcal{O}_N(f_N(\nu))\sum_{m=1}^{\infty}\eta_0^{(m+1)}+\mathcal{O}_N\left(\frac{1}{\nu^N}\right)\\
&=\sum_{m=1}^{\infty}\left[\left(\sum_{n=0}^{N-1}\sum_{k=0}^{n}\frac{A_{m+1,k}\left(\epsilon_1,\ldots,\epsilon_k\right)\eta_{n-k}^{(m+1)}}{\nu^n}\right)\right]+\mathcal{O}_N(f_N(\nu))+\mathcal{O}_N\left(\frac{1}{\nu^N}\right)\\
&=\sum_{n=0}^{N-1}\sum_{k=0}^{n}\left(\sum_{m=1}^{\infty}\frac{A_{m+1,k}\left(\epsilon_1,\ldots,\epsilon_k\right)\eta_{n-k}^{(m+1)}}{\nu^n}\right)+\mathcal{O}_N(f_N(\nu))+\mathcal{O}_N\left(\frac{1}{\nu^N}\right).
\end{align*}
Substituting these two expressions into \eqref{main_eqn}, we obtain that
\begin{align*}
&\frac{1}{2}=\frac{3c}{4}+\frac{3}{4\nu}\sum_{n=0}^{N-2}\left((-1)^{n+1}c+\sum_{m=0}^{n}(-1)^{n-m}\epsilon_{m+1}\right)\frac{1}{\nu^n}+\mathcal{O}_N(f_N(\nu))+\mathcal{O}_N\left(\frac{1}{\nu^N}\right)\\
&+\sum_{n=0}^{N-1}\sum_{k=0}^{n}\left(\sum_{m=1}^{\infty}\frac{A_{m+1,k}\left(\epsilon_1,\ldots,\epsilon_k\right)\eta_{n-k}^{(m+1)}}{\nu^n}\right)+\mathcal{O}_N(f_N(\nu))+\mathcal{O}_N\left(\frac{1}{\nu^N}\right)\\
&=\frac{3c}{4}+\frac{3}{4}\sum_{n=0}^{N-2}\left((-1)^{n+1}c+\sum_{m=0}^{n}(-1)^{n-m}\epsilon_{m+1}\right)\frac{1}{\nu^{n+1}}
+\sum_{m=1}^{\infty}A_{m+1,0}\eta_{0}^{(m+1)}\\
&\qquad\qquad+\sum_{n=1}^{N-1}\sum_{k=0}^{n}\left(\sum_{m=1}^{\infty}\frac{A_{m+1,k}\left(\epsilon_1,\ldots,\epsilon_k\right)\eta_{n-k}^{(m+1)}}{\nu^n}\right)+\mathcal{O}_N(f_N(\nu))+\mathcal{O}_N\left(\frac{1}{\nu^N}\right)\\
&=\frac{3c}{4}+\frac{3}{4}\sum_{n=0}^{N-2}\left((-1)^{n+1}c+\sum_{m=0}^{n}(-1)^{n-m}\epsilon_{m+1}\right)\frac{1}{\nu^{n+1}}
+\sum_{m=1}^{\infty}A_{m+1,0}\eta_{0}^{(m+1)}\\
&\qquad\qquad+\sum_{n=0}^{N-2}\sum_{k=0}^{n+1}\left(\sum_{m=1}^{\infty}\frac{A_{m+1,k}\left(\epsilon_1,\ldots,\epsilon_k\right)\eta_{n-k+1}^{(m+1)}}{\nu^{n+1}}\right)+\mathcal{O}_N(f_N(\nu))+\mathcal{O}_N\left(\frac{1}{\nu^N}\right)\\
&=\frac{3c}{4}+\sum_{m=1}^{\infty}A_{m+1,0}\eta_{0}^{(m+1)}+R_{N}(\nu)\sum_{m=1}^{\infty}\eta_0^{(m+1)}+\mathcal{O}_N(f_N(\nu))+\mathcal{O}_N\left(\frac{1}{\nu^N}\right)\\
&\qquad\qquad+\sum_{n=0}^{N-2}\left[\frac{3}{4}\left((-1)^{n+1}c+\sum_{m=0}^{n}(-1)^{n-m}\epsilon_{m+1}\right)+\sum_{k=0}^{n+1}\left(\sum_{m=1}^{\infty}A_{m+1,k}\left(\epsilon_1,\ldots,\epsilon_k\right)\eta_{n-k+1}^{(m+1)}\right)\right]\frac{1}{\nu^{n+1}}
\end{align*}
In view of the equations \eqref{epsilon_rln} and \eqref{c_value}, the above equation reduces to
\begin{align*}
	\mathcal{O}_N\left(\frac{1}{\nu^N}\right)+\mathcal{O}_N(f_N(\nu))=0.
%	\nu R_N(\nu)\eta_1(\nu)+\mathcal{O}_N\left(\frac{1}{\nu^N}\right)+R_{N}(\nu)\sum_{m=1}^{\infty}\eta_0^{(m+1)}=0
\end{align*}
Since $R_N(\nu)=\mathcal{O}_N(f_N(\nu))$ we conclude that
\begin{align*}
R_N(\nu)=\mathcal{O}_N\left(\frac{1}{\nu^N}\right)
\end{align*}
as $\nu\to \infty$.
Moreover, by substituting $n=0$ in \eqref{epsilon_rln} we obtain that
\begin{align*}
	 \frac{3}{4}\left(-c+\epsilon_1\right)+\sum_{m=1}^{\infty}\left(A_{m+1,0}\eta_1^{(m+1)}+A_{m+1,1}(\epsilon_1)\eta_0^{(m+1)}\right)=0.
\end{align*}
By using the fact that $A_{m+1,k}(\epsilon_1,\ldots,\epsilon_k)$ is the ordinary potential polynomial of $\left(c+\sum_{n=1}^{N-1}\frac{\epsilon_n}{\nu^n}\right)^{m+1},$ we
conclude \eqref{eps_1_value}. Similarly, we obtain the value of $\epsilon_i$ for $i \in \{2,3,\ldots\}$. This completes the proof of the
theorem.
\end{proof}	

\begin{proof}[\bf Proof of Theorem \ref{Theorem2}]
From the proof of \cite[Lemma 2.5]{BS14} we arrive at
\begin{equation}\label{h_exp} 1+\frac{zh_\nu^{\prime\prime}(z)}{h_\nu^\prime(z)}=1+\frac{z^{\frac{1}{2}}}{2}\frac{z^{\frac{1}{2}}J_{\nu+2}(z^{\frac{1}{2}})-4J_{\nu+1}(z^{\frac{1}{2}})}{2J_\nu(z^{\frac{1}{2}})-
z^{\frac{1}{2}}J_{\nu+1}(z^{\frac{1}{2}})}=1-\sum_{n=1}^{\infty}\frac{z}{\beta_{\nu,n}^2-z},
\end{equation}
where $\beta_{\nu,n}$ is the $n$th positive zero of the Dini function $e_\nu(z)$, defined in \eqref{dini_e}. For $\alpha=0$, \cite[Theorem 1.3]{BS14} implies that $z=r^c(h_\nu)$ is the smallest positive root of $$z\mapsto 1+\frac{z^{\frac{1}{2}}}{2}\frac{z^{\frac{1}{2}}J_{\nu+2}(z^{\frac{1}{2}})-4J_{\nu+1}(z^{\frac{1}{2}})}{2J_\nu(z^{\frac{1}{2}})-z^{\frac{1}{2}}J_{\nu+1}(z^{\frac{1}{2}})}.$$
Since the expression in \eqref{h_exp} is equal to zero at $z=r^c(h_\nu)$ we obtain that
\begin{equation}\label{eqn_beta}
	1-\sum_{n=1}^{\infty}\frac{r^c(h_\nu)}{\beta_{\nu,n}^2-r^c(h_\nu)}=0.
\end{equation}
Now with the help of Lemma \ref{Lemma4} we arrive at
$$r^c(h_\nu)=\nu\left(d+\mathcal{O}\left(\frac{1}{\nu}\right)\right)=\nu\left(d+\epsilon(\nu)\right)$$
for large $\nu$, where $d$ is some constant and $\epsilon(\nu)=\mathcal{O}\left(\frac{1}{\nu}\right)$. Rearranging \eqref{eqn_beta} we find that
\begin{align*}
1&=\sum_{n\geq1}^{}\frac{r^c(h_\nu)}{\beta_{\nu,n}^2-r^c(h_\nu)}=\sum_{n\geq1}^{}\frac{(\nu\left(d+\epsilon(\nu)\right))}{\beta_{\nu,n}^2-(\nu\left(d+\epsilon(\nu)\right))}\\
&=\sum_{n\geq1}^{}\frac{1}{\beta_{\nu,n}^2}\frac{\nu\left(d+\epsilon(\nu)\right)}{1-\frac{\nu\left(d+\epsilon(\nu)\right)}{\beta_{\nu,n}^2}}\\
&=\sum_{n\geq1}^{}\frac{\nu\left(d+\epsilon(\nu)\right)}{\beta_{\nu,n}^2}\sum_{m\geq0}^{}\frac{\left(\nu\left(d+\epsilon(\nu)\right)\right)^m}{\beta_{\nu,n}^{2m}}\\
&=\sum_{m\geq0}^{}\left(\nu\left(d+\epsilon(\nu)\right)\right)^{m+1}\sum_{n\geq1}^{}\frac{1}{\beta_{\nu,n}^{2m+2}}\\
&=\sum_{m\geq0}^{}\left(\nu\left(d+\epsilon(\nu)\right)\right)^{m+1}\theta_{m+1}(\nu),
\end{align*}
provided $\nu$ is sufficiently large. The rest of the proof is very similar to the proof of Theorem \ref{Theorem1} and hence we omit the details.
\end{proof}	

\begin{proof}[\bf Proof of Theorem \ref{unif_Theorem1}]
	From \cite[Theorem 1]{BPS14} for $\nu>-1$ the Dini function $d_\nu$, as defined in \eqref{dini_d}, admits the Weistrassian factorization as
	\begin{equation}\label{unif_d_nu_prod}
		d_\nu(z)=\frac{z^\nu}{2^\nu\Gamma\left(\nu+1\right)}\prod_{n=1}^{\infty}\left(1-\frac{z^2}{\alpha_{\nu,n}^2}\right),
	\end{equation}
	where $\alpha_{\nu,n}$ is the $n$th positive zeros of the Dini function $d_\nu(z)$ and the infinite product is uniformly convergent on each compact subset of the complex plane.	
	From \eqref{unif_g_nu_def} and \eqref{dini_d} we obtain that
	$$
	g_\nu^\prime(z)=2^\nu\Gamma(\nu+1)z^{-\nu}d_\nu(z).
	$$
	By using the above relation, the infinite product expression  \eqref{unif_d_nu_prod} and the expression \eqref{unif_g_nu_expres} we have
	\begin{align*}
		\frac{zg_\nu^{\prime\prime}(r)}{g_\nu^\prime(r)}=-\sum_{n=1}^{\infty}\frac{2z^2}{\alpha_{\nu,n}^2-z^2}=z\frac{(2\nu -1)J_{\nu+1}(z)-zJ_\nu(z)}{J_\nu(z)-zJ_{\nu+1}(z)}.
	\end{align*}
	By using \cite[eqn. 3.9]{DS17} and above relation we obtain
	\begin{align*}
		\real\left(1+\frac{zg_\nu^{\prime\prime}(z)}{g_\nu^\prime(z)}\right)-\left|\frac{zg_\nu^{\prime\prime}(z)}{g_\nu^\prime(z)}\right|&\geq 1+2\frac{rg_\nu^{\prime\prime}(r)}{g_\nu^\prime(r)}\\
		&=1+2r\frac{(2\nu-1)J_{\nu+1}(r)-rJ_{\nu}(r)}{J_\nu(r)-rJ_{\nu+1}(r)}\\
		%	&=1+r\frac{rJ_{\nu+2}(r)-3J_{\nu+1}(r)}{J_\nu(r)-rJ_{\nu+1}(r)}\\
		&=1-2\sum_{n=1}^{\infty}\frac{2r^2}{\alpha_{\nu,n}^2-r^2},
	\end{align*}
	for $|z|\leq r <\alpha_{\nu,1}$. Also from \cite[Theorem 3.2]{DS17} the radius of uniform convexity of $g_\nu(z)$ is the smallest positive root of the equation $$1+2r\frac{(2\nu-1)J_{\nu+1}(r)-rJ_{\nu}(r)}{J_\nu(r)-rJ_{\nu+1}(r)}=0.$$ From the above discussion and the definition of the radius of uniform convexity, at $r=r^{uc}(g_\nu)$ we obtain that
	\begin{equation}\label{unif_eqn24}
		1-4\sum_{n=1}^{\infty}\frac{(r^{uc}(g_\nu))^2}{\alpha_{\nu,n}^2-(r^{uc}(g_\nu))^2}=0.
	\end{equation}	
	Now with the help of Lemma \ref{unif_Lemma_g_asy} we have
	$$(r^{uc}(g_\nu))^2=\nu\left(\tilde{c}+\mathcal{O}\left(\frac{1}{\nu}\right)\right)=\nu\left(\tilde{c}+\varepsilon(\nu)\right)$$
	for large $\nu$, where $\tilde{c}$ is some constant and $\varepsilon(\nu)=\mathcal{O}\left(\frac{1}{\nu}\right)$. Rearranging \eqref{unif_eqn24} we find that
	\begin{align*}
		1&=4\sum_{n\geq1}^{}\frac{(r^{uc}(g_\nu))^2}{\alpha_{\nu,n}^2-(r^{uc}(g_\nu))^2}=4\sum_{n\geq1}^{}\frac{(\nu\left(\tilde{c}+\varepsilon(\nu)\right))}{\alpha_{\nu,n}^2-(\nu\left(\tilde{c}+\varepsilon(\nu)\right))}\\
		&=4\sum_{n\geq1}^{}\frac{1}{\alpha_{\nu,n}^2}\frac{\nu\left(c+\varepsilon(\nu)\right)}{1-\frac{\left(\tilde{c}+\varepsilon(\nu)\right)}{\alpha_{\nu,n}^2}}\\
		&=4\sum_{n\geq1}^{}\frac{\nu\left(\tilde{c}+\varepsilon(\nu)\right)}{\alpha_{\nu,n}^2}\sum_{m\geq0}^{}\frac{\left(\nu\left(\tilde{c}+\varepsilon(\nu)\right)\right)^m}{\alpha_{\nu,n}^{2m}}\\
		&=4\sum_{m\geq0}^{}\left(\nu\left(\tilde{c}+\varepsilon(\nu)\right)\right)^{m+1}\sum_{n\geq1}^{}\frac{1}{\alpha_{\nu,n}^{2m+2}}\\
		&=4\sum_{m\geq0}^{}\left(\nu\left(\tilde{c}+\varepsilon(\nu)\right)\right)^{m+1}\eta_{m+1}(\nu),
	\end{align*}
	which can be rewritten as
	\begin{equation}\label{unif_main_eqn}
		\frac{1}{4}=\nu\left(\tilde{c}+\varepsilon(\nu)\right)\eta_{1}(\nu)+\sum_{m=1}^{\infty}\left(\nu\left(\tilde{c}+\varepsilon(\nu)\right)\right)^{m+1}\eta_{m+1}(\nu),
	\end{equation}
	provided $\nu$ is sufficiently large. We would like to point out that \eqref{unif_main_eqn} is very similar to \eqref{main_eqn}. Consequently, the remainder of the proof closely parallels the proof of Theorem \ref{Theorem1}, and we therefore omit the details.
\end{proof}	

\begin{proof}[\bf Proof of Theorem \ref{unif_Theorem2}]
	
	From \cite[Lemma 2.5]{BS14} and \eqref{unif_h_nu_expres} we have
	\begin{align*}
		\frac{zh_\nu^{\prime\prime}(r)}{h_\nu^\prime(r)}=-\sum_{n=1}^{\infty}\frac{z}{\beta_{\nu,n}^2-z}=\frac{1}{2}z^{\frac{1}{2}}\frac{2(\nu -1)J_{\nu+1}(z^{\frac{1}{2}})-z^{\frac{1}{2}}J_\nu(z^{\frac{1}{2}})}{2J_\nu(z^{\frac{1}{2}})-z^{\frac{1}{2}}J_{\nu+1}(z^{\frac{1}{2}})},
	\end{align*}
	where $\beta_{\nu,n}$ is the $n$th positive zero of the Dini function $e_\nu(z).$
	By using \cite[eqn. (3.16)]{DS17} and the above relation we obtain that
	\begin{align*}
		\real\left(1+\frac{zh_\nu^{\prime\prime}(z)}{h_\nu^\prime(z)}\right)-\left|\frac{zh_\nu^{\prime\prime}(z)}{h_\nu^\prime(z)}\right|&\geq 1+2\frac{rh_\nu^{\prime\prime}(r)}{h_\nu^\prime(r)}\\
		&=1+r^{\frac{1}{2}}\frac{r^{\frac{1}{2}}J_{\nu+2}(r^{\frac{1}{2}})-4J_{\nu+1}(r^{\frac{1}{2}})}{2J_\nu(r^{\frac{1}{2}})-r^{\frac{1}{2}}J_{\nu+1}(r^{\frac{1}{2}})}\\
		&=1-\sum_{n=1}^{\infty}\frac{2r}{\beta_{\nu,n}^2-r},
	\end{align*}
	for $|z|<r<\beta_{\nu,1}^2$. We know from \cite[Theorem 3.3]{DS17} that the radius of uniform convexity of $h_\nu(z)$ is the smallest positive root of the equation
	$$1+r^{\frac{1}{2}}\frac{r^{\frac{1}{2}}J_{\nu+2}(r^{\frac{1}{2}})-4J_{\nu+1}(r^{\frac{1}{2}})}{2J_\nu(r^{\frac{1}{2}})-r^{\frac{1}{2}}J_{\nu+1}(r^{\frac{1}{2}})}=0.$$
	From the above discussion and the definition of the radius of uniform convexity, at $r=r^{uc}(h_\nu)$ we obtain that
	\begin{equation}\label{unif_eqn24_h}
		1-2\sum_{n=1}^{\infty}\frac{r^{uc}(h_\nu)}{\beta_{\nu,n}^2-r^{uc}(h_\nu)}=0.
	\end{equation}	
	Now, with the help of Lemma \ref{unif_Lemma_h_asy} we have
	$$r^{uc}(h_\nu)=\nu\left(\tilde{d}+\mathcal{O}\left(\frac{1}{\nu}\right)\right)=\nu\left(\tilde{d}+\varepsilon(\nu)\right)$$
	for large $\nu$, where $\tilde{d}$ is some constant and $\varepsilon(\nu)=\mathcal{O}\left(\frac{1}{\nu}\right)$. Rearranging \eqref{unif_eqn24_h} we find that
	\begin{align*}
		\frac{1}{2}&=\sum_{n\geq1}^{}\frac{r^{uc}(h_\nu)}{\beta_{\nu,n}^2-r^{uc}(h_\nu)}=\sum_{n\geq1}^{}\frac{(\nu\left(\tilde{d}+\varepsilon(\nu)\right))}{\beta_{\nu,n}^2-(\nu\left(\tilde{d}+\varepsilon(\nu)\right))}\\
		&=\sum_{n\geq1}^{}\frac{1}{\beta_{\nu,n}^2}\frac{\nu\left(\tilde{d}+\varepsilon(\nu)\right)}{1-\frac{\nu\left(\tilde{d}+\varepsilon(\nu)\right)}{\beta_{\nu,n}^2}}\\
		&=\sum_{n\geq1}^{}\frac{\nu\left(\tilde{d}+\varepsilon(\nu)\right)}{\beta_{\nu,n}^2}\sum_{m\geq0}^{}\frac{\left(\nu\left(\tilde{d}+\varepsilon(\nu)\right)\right)^m}{\beta_{\nu,n}^{2m}}\\
		&=\sum_{m\geq0}^{}\left(\nu\left(\tilde{d}+\varepsilon(\nu)\right)\right)^{m+1}\sum_{n\geq1}^{}\frac{1}{\beta_{\nu,n}^{2m+2}}\\
		&=\sum_{m\geq0}^{}\left(\nu\left(\tilde{d}+\varepsilon(\nu)\right)\right)^{m+1}\theta_{m+1}(\nu),
	\end{align*}
	provided $\nu$ is sufficiently large. The rest of the proof of is very similar to the proof of Theorem \ref{unif_Theorem1} and hence we omit the details.
\end{proof}	
\subsection*{Acknowledgments} Pranav Kumar is grateful to the Council of Scientific and Industrial Research India (Grant No. 09/1022(0060)/2018-EMR-I) for the financial support and Sanjeev Singh is thankful to the Science and Engineering Research Board (SERB), Department of Science and Technology, Government of India for the financial support through Project MTR/2022/000792.

\end{document}